\documentclass[10pt]{IEEEtran}

\usepackage{verbatim}
\usepackage{amsfonts,amsmath,mathrsfs,amssymb,amsbsy}
\usepackage[final]{graphicx}
\usepackage{times,cite}

\usepackage{balance}

\long\def\symbolfootnote[#1]#2{\begingroup%
\def\thefootnote{\fnsymbol{footnote}}\footnote[#1]{#2}\endgroup}

\usepackage{verbatim}
\usepackage[]{float,latexsym}
\usepackage{amsfonts,amsmath,mathrsfs,amssymb,amsbsy}
\usepackage{url}

\newtheorem{theorem}{Theorem}[section]

\newtheorem{lemma}{Lemma}

\newtheorem{assumption}{Assumption}

\newcommand{\Prob}{\mathsf{P}}
\newcommand{\Expect}{\mathsf{E}}

\usepackage{color}
\definecolor{lightblue}{rgb}{.7, .8, 1}
\definecolor{lightgreen}{rgb}{.6, 1, .6}
\usepackage{color}
\definecolor{brown}{rgb}{1,0.38,0.03}

\definecolor{OliveGreen}{rgb}{.2,0.6,0.2}
\definecolor{BrickRed}{rgb}{.7,0.2,0.2}

\newcommand{\ignore}[1]{} 

\long\def\symbolfootnote[#1]#2{\begingroup%
\def\thefootnote{\fnsymbol{footnote}}\footnote[#1]{#2}\endgroup}

\DeclareMathOperator*{\argmax}{arg\,max}

\newcommand{\taug}{\tau_{\scriptscriptstyle \text{G}}}

\def\bR{{\mathbf{R}}}

\def\bR{{\mathbf{R}}}

\def\bX{\mathbf{X}}
\def\bZ{\mathbf{Z}}

\def\bmu{{\mathbf{\mu}}}

\begin{document}

\sloppy

\title{Quickest Detection for Changes in Maximal kNN Coherence of Random Matrices}

\author{Taposh Banerjee, Hamed Firouzi, and Alfred O. Hero III, \textit{Fellow, IEEE \vspace{-3ex}}
}

\maketitle
\symbolfootnote[0]{\small
This work was partially supported by the Consortium for Verification Technology under Department of Energy National Nuclear Security Administration award number DE-NA0002534 and the Air Force Office of Scientific Research under grant FA9550-13-1-0043.
A preliminary version of this paper has been presented at 2015 IEEE International Symposium
on Information Theory.

Taposh Banerjee is with the School of Engineering and Applied Sciences, Harvard University,
		Cambridge, MA, Hamid Firouzi is with Goldman Sachs, New York, NY, and Alfred Hero is with
		the Department of Electrical Engineering and Computer Science, University of Michigan, Ann Arbor, MI. This work was completed when the authors were with the Department of Electrical Engineering and Computer Science, University of Michigan, Ann Arbor, MI 48109, USA
(e-mail: tbanerjee@seas.harvard.edu; firouzi@umich.edu; hero@umich.edu)

}

\begin{abstract}
This paper addresses the problem of quickest detection of a change in the maximal coherence between columns of a $n\times p$ random matrix based on a sequence of matrix observations having a single unknown change point. 
The random matrix is assumed to have identically distributed rows and the  maximal coherence is defined as the largest of the $p \choose 2$ correlation coefficients associated with any row.   Likewise the $k$ nearest neighbor (kNN) coherence is defined as the $k$-th largest of these correlation coefficients.
The forms of the pre- and post-change distributions of the observed matrices are assumed to belong to the family of elliptically contoured densities with sparse dispersion matrices but are otherwise unknown.  A non-parametric stopping rule is proposed that is based on the maximal k-nearest neighbor sample coherence between columns of each observed random matrix. This is a summary statistic that is related to a test of existence of a hub vertex in a sample correlation graph having degree at least $k$. Performance bounds on the delay and false alarm performance of the proposed stopping rule are obtained in the purely high dimensional regime where $p\rightarrow \infty$ and $n$ is fixed.
When the pre-change dispersion matrix is diagonal it is shown that, among all functions of the proposed summary statistic, the proposed stopping rule is asymptotically optimal under a minimax quickest change detection (QCD) model as the stopping threshold approaches infinity. 
The theory developed also applies to sequential hypothesis testing and fixed sample size tests.
\end{abstract}



\section{Introduction} \label{sec:Intro}
One of the greatest challenges in data analysis is to develop robust algorithms for statistical inference on large scale data.
Many big data applications fall in the so-called sample starved regime \cite{hero-bala-submitted-2015}, where conclusions have to be drawn or decisions
have to be made based on a small set of samples of a high-dimensional vector.
Most classical statistical tests have been designed for the large sample regime,
where the number of samples are much larger than the dimension of the vector,
and hence are not applicable to
high-dimensional data analysis. Thus, new approaches are needed to address these challenges.

In this paper we consider the problem of detecting a change in maximal kNN coherence
between 
columns of a random  matrix. The kNN coherence is defined as the $k$-th largest correlation coefficient between columns of the random matrix, as described below, and arises in dependency testing ($k=1$), correlation screening ($k=1$), node centrality analysis and hub discovery ($k>1$) in graphical models. For simplicity, unless it might cause confusion, we refer to maximal kNN coherence as simply maximal coherence. 

We take a non-parametric approach in this paper, assuming that the data matrix has a distribution in the elliptically contoured family but is otherwise unknown.  This non-parametric family is very general, containing both the light-tailed matrix normal (Gaussian) distribution and the heavy tailed multivariate Student-$t$ distribution.   The elliptically contoured  family of distributions is characterized by a vector valued location parameter, a vector valued scale parameter, a matrix valued  dispersion parameter, and a shaping function. The coherence matrix is derived from the dispersion matrix by dividing the $i$-th and $j$-th row and column by the square root of the product of the $i$th and $j$-th diagonal element, resulting in a matrix of correlation coefficients.  We will use the sample correlation coefficient, computed from pairs of columns of the data matrix, as a change detection statistic. Under the elliptically contoured assumption on the observed sequence of  matrices, this empirical estimate of coherence has a  distribution that only depends on the dispersion parameter; it does not depend on the location parameter nor the shaping parameter.

In the maximum coherence quickest change detection setting addressed here, the first few matrix-valued observations in the sequence are  
i.i.d. with a nominal coherence matrix, i.e., a ``normal'' or ``expected'' baseline of multivariate correlations.
At some time point in the sequence the maximum coherence may change, e.g., a sudden shift from incoherence to  coherence, or vice-versa.
The objective is to detect a change in maximal coherence as quickly as possible.
In many applications the change has to be detected in real time, i.e., with minimum possible delay, 
while avoiding false alarms.
Rapid and timely detection of
disorder can potentially save the cost of acquiring the rest of the samples.

The maximal coherence change detection problem has many applications including to slippage problems in multivariate time-series analysis and financial stock analysis, anomaly detection in social networks, cyber-physical systems and communication networks, and intrusion detection in sensor networks.
In multivariate time-series analysis, it is of interest to know if the maximal correlations between time series have abruptly changed over time. 
In portfolio balancing, it is of interest
to detect a sudden change in the maximal  correlation between a set of stocks being monitored.
In social networks, it is of interest to detect an abrupt change in the interaction levels between agents.
In communication networks it is of interest to detect emergent hubs of highly correlated traffic flows over
the network, which may be a potential point of attack by a cyber attacker. In sensor network intrusion detection,
the presence of an intruder may abruptly increase or decrease maximal correlation between various sensors located near the intruder.

The major challenges in this problem are:
\begin{enumerate}
\item In the high dimensional setting the number $p$ of variables (columns) in the data matrix may vastly exceed the number $n$ of samples (rows) in the matrix.
\item The statistical properties of each data matrix may not be precisely known, i.e., the problem is nonparametric in nature\footnote{By nonparametric we mean that the parameter space is infinite dimensional.}.
\item Accurate detection of a change in the maximal correlation is complicated by the unknown baseline (pre-change) distribution of data matrix, and in particular the maximal sample correlation test statistic.    
\end{enumerate}
We overcome these challenges by using recent results in random matrix theory to obtain an asymptotic large $p$ small $n$ distribution of the maximal sample correlation and recent results in sequential detection to define change detectors whose false alarm rates can be controlled even when the nominal (pre-change) distribution has unknown parameters.

We formulate the change detection problem in the framework of quickest change detection (QCD)
(see, e.g., \cite{veer-bane-elsevierbook-2013}, \cite{poor-hadj-qcd-book-2009},
and \cite{tart-niki-bass-2014}). In the QCD problem a decision maker observes a stochastic process over time.
At some point in time, called the change point, the distribution of the
process changes. The decision maker has to detect this change in distribution with
minimum possible delay, subject to a constraint
on false alarms. The QCD problem has been formulated in various ways in the literature.
One prevalent formulation of the QCD problem is as a stochastic optimization problem,
where the goal is to find a stopping time
on the observed stochastic process so as to minimize a suitable metric on the delay, subject to a suitable
metric on the false alarm rate. A typical solution is a stopping rule that reduces to a single threshold test,
where a sequence of statistics is computed over time, and a change is declared the first time the statistic exceeds a stopping threshold.
The stopping threshold is chosen to control the rate of false alarms.
The theoretical foundation for such sequential decision making was laid by Wald; see \cite{wald_book_1947}
\cite{wald-wolf-amstat-1948}.
The Bayesian version of the problem, where a prior on the change point is assumed, is developed
in \cite{girs_rubin_1952}, \cite{shir-siamtpa-1963}, and \cite{tart-veer-siamtpa-2005}.
The QCD problem in non-Bayesian minimax settings has been formulated in \cite{page-biometrica-1954}, \cite{lord-amstat-1971},
\cite{mous-astat-1986}, \cite{ritov-astat-1990}, \cite{lai-ieeetit-1998}.
In general, an optimal or asymptotically optimal solution to a QCD problem can be obtained only when the
pre- and post-change distributions are known to the decision maker, or when
the post-change distribution is in a parametric family.
In the nonparametric setting, an optimal solution is hard to obtain.
As a result, in the nonparametric setting considered here  the goal is often
less ambitious than to find an optimal solution. Rather, a reasonable procedure is proposed and its properties are established, e.g., consistency, convergence rate, scalability, and so on.
In this paper we propose a consistent and scalable nonparametric procedure for correlation change detection in 
a high-dimensional sample starved setting.
See Section~\ref{sec:Prob} and \ref{sec:SumStat} 
for details. 

Specifically, we consider the following random matrix observation model.
An independent sequence of random matrices $\{\mathbb{X}(m)\}$ is observed over time, indexed by $m$,
where each $\mathbb{X}(m)$ is an $n \times p$ short and fat random matrix.
By short and fat matrix
we mean $p \gg n$.
The $n$ rows of each random matrix may correspond to a block of identically distributed random samples of a $p$ variate vector, e.g.,  sampled over blocks of time or sampled in a sequence of repeated experiments.
For example, in the case of detecting a change in the coefficients of a Gaussian univariate time series,
$p$ successive time samples may be acquired over $n$ well separated blocks of time. 
In multi-pulse radar imaging used for change detection, $p$ denotes the number of pixels and $n$ denotes the number of radar pulses used in a single radar burst to form the radar image at time $m$.
In stochastic finance, we may have access to multiple instances of stock returns over a week.
A quarterly report consists of a matrix of $n= 12$ weekly returns of the $p=10,000$ stocks. The objective is to detect a change in the maximum correlation associated with the sequence of quarterly reports. 
In each of these examples the data has been grouped into batches of time to form the data matrix $\mathbb{X}$. Our analysis in this paper will require that the change in distribution occur between batches rather than within a batch.

If the distribution of the random matrices belong to a parametric family,
and the value of the parameter before the change is known, then efficient procedures
from the quickest change detection literature can be used for detection \cite{lord-amstat-1971}, \cite{lai-ieeetit-1998}, \cite{lord-poll-nonantiest-2005}.
However, as discussed above, 
 when the pre-change parameter or the distribution is unknown,
no optimal procedures are known for detection of change. Here, by optimal we mean optimal in the sense of minimizing
detection delay as studied in the classical QCD literature \cite{shir-siamtpa-1963}, \cite{lord-amstat-1971}, \cite{poll-astat-1985}. 
In this paper we propose a nonparametric procedure that can provably detect a 
 a change in the maximal coherence in the high dimensional regime of $p\gg n$ and $n$ fixed.
 
Specifically, in Section~\ref{sec:Prob} we consider the problem of quickest detection of a change in
maximal coherence
under the assumption
that the $\mathbb{X}(m)$ are independent and identically distributed,
with joint distribution from the nonparametric family of elliptically contoured distributions.
We propose the maximal sample correlation statistic $V(\mathbb{X})$
as the test statistic derived from the data matrix $\mathbb{X}$ that is used to detect a change.
The summary statistic 
$V(\mathbb X)$ is an estimate of the $k$-th largest correlation coefficient in the ensemble coherence matrix. While in the classical setting of $p$ fixed and $n\rightarrow \infty$ this estimate is asymptotically consistent, it has no theoretical accuracy guarantees in our high dimensional setting of $p\rightarrow \infty$ and $n$ fixed. However, by interpreting $V(\mathbb X)$ as the  
minimal size of the $k$-nearest neighborhood among all
the columns of the observed matrix, results of \cite{hero-bala-IT-2012} establish an explicit parametric  form for the distribution of $V(\mathbb X)$ in the high dimensional setting that is the basis for our proposed QCD procedure. 
Here the size of the neighborhood is measured by the sample correlation associated with the column and its $k$-nearest (most correlated) neighbors.
In particular, we obtain an asymptotic distribution for the summary statistic in
the sample starved purely high dimensional regime of $p \rightarrow \infty$ with $n$ fixed and small.
We show in Section~\ref{sec:SumStat} that the distribution of the summary statistic belongs to a one-parameter exponential family,
with the unknown parameter a function of the underlying distribution of the sample coherence matrix.

In this manner we map the sequence of $n\times p$ observed data matrices $\{\mathbb{X}(m)\}$ to a
sequence of real valued summary statistics  $\{V(\mathbb{X}(m))\}$ whose distribution is in a known parametric family for sufficiently large $p$ and finite $n$.
Any change in distribution in the sequence $\{\mathbb{X}(m)\}$ that causes a change in the parameter of the distribution
of  $\{V(\mathbb{X}(m))\}$ can be detected. 
While the parametric family for $\{V(\mathbb{X}(m))\}$ is known, the actual pre- and post-change parameters may not be known.
However, under block sparsity conditions on 
the population correlation matrix, the pre-change parameter can be shown to be in the vicinity of $1$
(See Theorem~\ref{thm:CorrScr} below). 
We thus detect the change in parameter by applying the generalized likelihood ratio (GLR) based test of 
Lorden \cite{lord-amstat-1971}  (also see \cite{lai-jrss-1995})
by setting the pre-change parameter to $1$. 
We then provide a detailed performance analysis of the Lorden's test under misspecification of the pre- and post-change 
distributions. 
We emphasize that the analysis for Lorden's GLR test is non-trivial and requires additional conditions beyond the standard conditions used for SPRT and CUSUM \cite{tart-niki-bass-2014}
(see Section~\ref{sec:unknownJ0}).

We show that the proposed test procedure is asymptotically optimal over all tests that use the summary statistic  $V(\mathbb X)$ to detect changes in maximal correlation. 
 This is because our proposed test uses the asymptotic distribution of the summary statistic in the Lorden's test.

We also remark that while we focus on sequential detection of changes, the proposed summary statistic can also 
be used for more general sequential testing problems in the high dimensional data regime \cite{tart-niki-bass-2014}.

In Section~\ref{sec:Numerical} we validate the effectiveness of the proposed procedure by verifying the theoretical results through numerical simulations.

In summary, the contributions of this paper: we propose a nonparametric quickest detection procedure for detecting a change in correlation in a sequence of $n \times p$ high-dimensional matrices ($p \gg n$) $\{\mathbb{X}(m)\}$. Specifically, the contributions of this paper are as follows:
\begin{enumerate}
  \item We propose a novel summary statistic $V(\mathbb{X})$ of the data matrix $\mathbb{X}$ to test maximal coherency. This statistic is the maximal sample correlation between the columns of the matrix $\mathbb{X}$.
  \item We obtain an asymptotic distribution for $V(\mathbb{X})$ in the purely high dimensional regime of fixed $n$ and 
   $p \to \infty$ under the assumption that the matrices $\mathbb X$ are elliptically distributed. This asymptotic distribution of $V(\mathbb{X})$ belongs to a one-parameter exponential family.
   \item We propose applying Lorden's GLR test to the summary statistic sequence $\{V(\mathbb{X}(m))\}$ and we
analyze its performance when the pre- and post-change distribution are misspecified. 
  \item We obtain conditions on the  pre- and post-change dispersion matrices in the elliptical model of $\mathbb X(m)$ for which the change can be accurately detected for any specified false positive constraint.
\end{enumerate}

\section{Problem Description}\label{sec:Prob}
A decision-maker sequentially acquires samples from a family of distributions of
$n \times p$ random matrices over time, indexed by $m$, leading to the random matrix sequence $\{\mathbb{X}(m)\}_{m \geq 1}$,
called data matrices.
The data matrices in the sequence are assumed to be statistically independent. 
For each $m$ the data matrix $\mathbb{X}(m)$ has the following properties. Each of its
$n$ rows is a sample
of a $p$-variate random vector $\mathbf{X}(m)$ $=[X_1(m), \cdots, X_p(m)]^T$ with $p \times 1$
location parameter $\boldsymbol{\mu}_m$ and
$p\times p$ positive definite dispersion matrix $\mathbf{\Sigma}_m$.
The random matrix $\mathbb X(m)$ is assumed to have a vector elliptically contoured distribution,
also called an elliptical density, \cite{ander-mutlstat-1996}, \cite{AndersonTR}, 
\begin{equation}\label{eq:matrixdensityexpr}
f_{\mathbb{X}(m)}(\mathbb{X}) = h_m(\text{trace}\{(\mathbb{X}-{\mathbf 1}{\mathbf \mu_m}^T )^T\boldsymbol{\Sigma}_m^{-1}(\mathbb{X}-{\mathbf 1}{\mathbf \mu_m}^T)\}),
\end{equation}
for shaping function $h_m$ on $\mathbb{R}^+$. 
The shaping function only has to satisfy the properties: 1) $h_m\geq 0$; 2) $\int h_m(u) du <\infty$ for \eqref{eq:matrixdensityexpr} to be a properly defined density function over $\mathbb R^p$.
If $\boldsymbol{\mu}_m=0$
and $\mathbf{\Sigma}_m = \mathbf{I}_p$, where $\mathbf{I}_p$ is the $p \times p$ identity matrix, then
the matrix $\mathbb X(m)$ is said to have a spherical density.
Note that under the density  \eqref{eq:matrixdensityexpr} the rows of $\mathbb X(m)$ are not independent unless the shaping function has the form $h_m(u)=c\exp(-u^2/2)$ for which $\mathbb X(m)$ is multivariate Gaussian.   When $h_m$ satisfies additional conditions the location parameter $\mathbf \mu_m$ and the dispersion parameter $\mathbf \Sigma_m$ are the mean and covariance of $f_{\mathbf X(m)}$. However, we do not require any additional conditions on $h_m$ in this paper.

As in \cite{hero-bala-IT-2012} and \cite{hero-bala-jasta-2011}, the $p\times p$ dispersion matrix $\mathbf \Sigma_m$ is said to be block sparse of degree $s$ if there exists a row-column permutation that puts it into block diagonal form with block size of size $s \times s$ where  $s=o(p)$.  It is said to be row-sparse of order $s$ if no row contains more than $s$ non-zero elements, where $s=o(p)$. Block sparsity, or the weaker row sparsity, of the dispersion matrix will be required for most of the results in this paper.

The coherency matrix, also known as the correlation matrix, is obtained from the dispersion matrix by normalizing the entries of $\mathbf \Sigma_m$ so that the diagonals are equal to 1. Specifically, the coherence matrix $\mathbf C_m$ is defined as:  $\mathbf C_m=\mathbf{D}_m^{-1} \mathbf{\Sigma}_m \mathbf{D}_m^{-1}$, where $\mathbf{D}_m$ is the diagonal matrix obtained by setting the off-diagonals of $\mathbf \Sigma_m$ equal to zero.  The magnitudes of the off-diagonal elements of $\mathbf C_m$ are the coherence coefficients of the model \eqref{eq:matrixdensityexpr}. 

The objective in this paper is to detect a change in the $k$-th largest coherence coefficient. For  $k=1$ this is simply the maximum  $\max_{i>j} {\mathrm {abs}}([[\mathbf C_m]]_{ij})$, where $[[\mathbf C_m]]_{ij}$ denotes its $i,j$th element, and abs$(u)$ denotes the absolute value of $u\in \mathbb R$.  As the maximum is equal to zero if and only if $\mathbf \Sigma_m$ is diagonal, it is relevant to multivariate dependency testing \cite{moran1980testing}, \cite{cameron1985new},\cite{anandkumar2009detection}. Similarly, as the maximal kNN coherence is equal to zero if and only if $\mathbf{\Sigma}_m$ is row sparse (with no row having more than $k$ non-zero entries) it is relevant to $k$-dependency testing or, equivalently, detection of high degree $(\geq k)$ vertices in a correlation graphical model \cite{cox2014multivariate},\cite{hero-bala-IT-2012}.  

With $\eta_m$ denoting the maximal kNN coherence of $\mathbb X(m)$  (here the dependence on $k$ is supressed), the objective of the change detection problem is to detect a change point $\gamma$ such that 
\begin{eqnarray}
\eta_m&=&\eta_0, \hspace{0.1in} \mathrm{for} \; m< \gamma 
\\
\eta_m&=&\eta_1, \hspace{0.1in} \mathrm{for} \; m\geq \gamma 
\end{eqnarray}
where $\eta_0\neq \eta_1$.

For the purpose of detecting and localizing a possible change in the maximal coherence, we will construct a test on a sequence of summary statistics  $\{V(m)\}_m$, where $V(m)$ is a function of the data matrix: $V(m)=V(\mathbb X(m))$. We elaborate on the form of the proposed summary statistic in Sec. III.  Let $f_V^0$ and $f_V^1$ denote the pre-change and post-change distributions of $V(m)$, i.e., $f_{V(m)}=f_V^0$ for $m<\gamma$ and $f_{V(m)}=f_V^1$ for $m\geq\gamma$.
At each time point $m$ the decision-maker decides to either stop sampling,
declaring that the change has occurred, i.e., $m\geq \gamma$, or to continue sampling.
The decision to stop at time $m$ is only a function of $(V(1), \cdots, V(m))$.
Thus, the time at which the decision-maker decides to stop sampling is a stopping time
for the sequence $\{V(m)\}$.
The decision-maker's objective is to detect the change in maximal correlation as quickly as possible,
subject to a constraint on the false alarm rate.

The above detection problem is an example of the quickest change detection (QCD) problem
\cite{poor-hadj-qcd-book-2009}, \cite{veer-bane-elsevierbook-2013}, and \cite{tart-niki-bass-2014}. 
The objective in our QCD problem is to find a stopping time $\tau$
on the sequence $\{V(m)\}$, so as to minimize a suitable metric on the delay $(\tau-\gamma)$,
subject to a constraint on a suitable metric on the event of false alarm $\{\tau < \gamma\}$.
This paper follows the minimax QCD formulation of Pollak \cite{poll-astat-1985}:
\begin{equation}\label{prob:Pollak}
\begin{split}
 \min_\tau       &  \quad \sup_{\gamma\geq 1}  \; \Expect_\gamma[\tau-\gamma| \tau \geq \gamma] \\
 \mbox{subj. to}  &  \quad \Expect_\infty[\tau] \geq \beta,
\end{split}
\end{equation}
where $\Expect_\gamma$ is the expectation with respect to the probability measure under which the change occurs at $\gamma$,
$\Expect_\infty$ is the corresponding expectation when the change never occurs, and $\beta \geq 1$ is a
user-specified constraint on the mean time to false alarm. 
Depending on the nature of the summary statistic $V(\mathbb X)$ the expectation $\Expect_\gamma$ may be 
a function of the shaping function, the translation parameters $\{\mu_m\}$, and the dispersion matrix $\{\Sigma_m\}$. 
For example, when the summary statistic is defined as the maximal sample coherency (see Section~\ref{sec:SumStat}), the expectation depends on the dispersion matrix through a scalar parameter $J$, but does not depend on the shaping function or the location vector. 
For simplicity of notation, we do not explicitly show the dependence of this expectation on these parameters.

If the pre- and post-change densities are known to the decision maker
then algorithms like the Cumulative Sum (CuSum) algorithm \cite{page-biometrica-1954}, \cite{lord-amstat-1971},
\cite{mous-astat-1986}, or the Shiryaev-Roberts family of algorithms \cite{robe-technometrics-1966}, \cite{poll-astat-1985},
\cite{tart-etal-thirdord-2012}, can be used to detect a change in the maximal coherence.. Both the CuSum algorithm
and the Shiryaev-Roberts family of algorithms have strong optimality properties with respect to both the popular
formulations of Lorden \cite{lord-amstat-1971} and that of Pollak \cite{poll-astat-1985}, used in this paper.

If the pre-change distribution is known and the post-change distributions belongs to a known parametric family, then
efficient QCD algorithms can be designed having strong asymptotic optimality properties,
based on, e.g., the generalized likelihood ratio (GLR) technique \cite{tart-niki-bass-2014},
the mixture based technique  \cite{tart-niki-bass-2014}, or
the non-anticipating estimation based technique \cite{lord-poll-nonantiest-2005}.


The summary statistic $V(m)=V(\mathbb X(m))$ that we will use for quickest change change detection is the $k$-th largest sample correlation, defined below, between the columns of the data matrix $\mathbb{X}(m)$. This is equivalent to the maximal kNN distance between columns, as measured by correlation distance.  
The theory from \cite{hero-bala-IT-2012} helps us establish that the proposed summary statistic
has a well defined exponential limiting distribution as $p\rightarrow \infty$ for fixed $n$, the so-called "purely high dimensional regime"
\cite{hero-bala-submitted-2015}.
This summary statistic is related to the empirical distribution of the vertex degree
of the correlation graph associated with the thresholded sample correlation matrix.
Below we show that the distribution of the statistic $V(\mathbb{X})$ converges to a
parametric distribution in the exponential family in this purely high dimensional regime.
Thus, the nonparametric QCD problem in terms of $\{\mathbb{X}(m)\}$ is mapped to a parametric QCD
problem in terms of the summary statistic sequence $\{V(\mathbb{X}(m))\}$.
We then apply a GLR based test suggested by Lorden in \cite{lord-amstat-1971} to the sequence of summary statistics
$\{V(\mathbb{X}(m))\}$ to detect the change efficiently.

If the pre-change dispersion matrix $\mathbf{\Sigma}_0$ is diagonal,
then we show below that the pre-change distribution $f_V^{0}$ is completely specified and the post-change distribution $f_V^{1}$ is specified up to a scalar unknown parameter.
In this case the GLR stopping rule used is asymptotically optimal under the
Lorden minimax QCD model \cite{lord-amstat-1971},
and hence also in terms of solving \eqref{prob:Pollak}, among all rules that
are stopping rules for the sequence $\{V(\mathbb{X}(m))\}$.

If the pre-change matrix $\mathbf{\Sigma}_0$ is unknown and not diagonal, then 
the pre-change distribution $f_V^{0}$ depends on an unknown scalar parameter.
Below we establish conditions on the matrix $\mathbf{\Sigma}_0$ which guarantee that
the GLR stopping rule remains approximately optimal, in the sense that the mean time to false alarm and mean time to detect are within a constant scaling factor of the values of the optimal QCD decision rule. This is achieved by analyzing the performance
of the GLR test under mis-specification of the pre-change distribution.

\section{Summary Statistic for Detecting a Change in Maximal Coherence}\label{sec:SumStat}
In this section we define the proposed summary statistic $V(\mathbb{X})$ and
obtain its asymptotic density in the
purely high dimensional regime of  $p\rightarrow \infty$, $n$ fixed. 

The notation below follows the conventions of \cite{hero-bala-IT-2012}.
For an elliptically distributed random data matrix $\mathbb{X}$ we write
$$\mathbb{X} = [\mathbf{X}_1, \cdots, \mathbf{X}_p] = [\mathbf{X}^T_{(1)}, \cdots, \mathbf{X}^T_{(n)}]^T,$$
where $\mathbf{X}_i = [X_{1i}, \cdots, X_{ni}]^T$ is the $i^{th}$ column and $\mathbf{X}_{(i)} = [X_{i1}, \cdots, X_{ip}]$
is the $i^{th}$ row.
Define the $p \times p$ sample covariance matrix as
$$ \mathbf{S} = \frac{1}{n-1} \sum_{i=1}^n (\mathbf{X}_{(i)} - \bar{\mathbf{X}})^T (\mathbf{X}_{(i)} - \bar{\mathbf{X}}),$$
where $\bar{\mathbf{X}}$ is the sample mean of the $n$ rows of $\mathbb{X}$.
Define the sample correlation (coherency) matrix as
$$ \mathbf{R} = \mathbf{D_S}^{-1/2} \mathbf{S} \mathbf{D_S}^{-1/2},  $$
where $\mathbf{D_A}$ denotes the matrix obtained by zeroing out all but the diagonal elements of the matrix $\mathbf{A}$.
Note that, under our assumption that the ensemble dispersion matrix $\mathbf{\Sigma}$ of the rows of $\mathbb X$ is positive definite, $\mathbf{D_S}$ is invertible with probability one.
Thus $\mathbf{R}_{ij}$, the $ij$th element of the matrix $\mathbf{R}$, is the
sample correlation coefficient between the $i^{th}$ and $j^{th}$ columns of $\mathbb X$.

Define $d^{(k)}_{\text{NN}}(i)$ to be the magnitude of the sample correlation between the $i$-th column of $\mathbb X$ and its
$k$-th nearest neighbor in the remaining columns of $\mathbb X$ (with respect to correlation distance). 
Equivalently, $d^{(k)}_{\text{NN}}(i)$ is defined in terms of the sample correlation::
$$d^{(k)}_{\text{NN}}(i) := k^{th} \mbox{ largest order statistic of } \{ |\mathbf{R}_{ij}|; j \neq i\}.$$
Then for fixed $k$, define the summary statistic
\begin{equation}\label{eq:Sum_Stat_V}
V_k(\mathbb{X}):= \max_{1 \leq i \leq p} d^{(k)}_{\text{NN}}(i).
\end{equation}
Thus, if magnitude correlation between variables is used as a distance measure, then 
the summary statistic $V_k$ is the maximum size of the k-Nearest Neighborhood (kNN) across variables. 
Note that the summary statistic $V_k$ is a global statistic, and is
insensitive to variations in the minimal kNN correlations as long as the maximum of these kNN correlations remains the same.

Below we show that the distribution of the statistic $V_k$ can be related to the distribution
of an integer valued random variable $N_{k, \rho}$ that counts the number of highly correlated neighborhoods.


For a threshold parameter $\rho \in [0,1]$ define the correlation graph $\mathcal{G}_\rho(\mathbf{R})$
associated with the correlation matrix $\mathbf{R}$
as an undirected graph with $p$ vertices, each representing a column of the data matrix $\mathbb{X}$.
An edge is present between vertices $i$ and $j$ if the magnitude of the sample correlation coefficient
between the $i^{th}$ and $j^{th}$ components of the random vector $\mathbf{X}$
is greater than $\rho$, i.e., if $|\mathbf{R}_{ij}| \geq \rho$, $i \neq j$.
We define $\delta_i$ to be the degree of vertex $i$ in the graph $\mathcal{G}_\rho(\mathbf{R})$.
For a positive integer $k \leq p-1$ we say that a vertex $i$ in the graph $\mathcal{G}_\rho(\mathbf{R})$
is a hub of degree $k$ if $\delta_i \geq k$. We denote
by $N_{k, \rho}$ the total number of hubs in the correlation graph $\mathcal{G}_\rho(\mathbf{R})$, i.e.,
$$N_{k, \rho} = \text{card} \{i: \delta_i \geq k \}. $$

The events $\{V_k(\mathbb{X}) \geq \rho\}$ and $\{N_{k, \rho} > 0\}$ are equivalent. Hence
\begin{equation}\label{eq:V_eq_N}
\Prob(V_k(\mathbb{X}) \geq \rho) = \Prob(N_{k, \rho} > 0).
\end{equation}

Because of the above relation, for a fixed level $\rho$,
$V_k(\mathbb{X})$ indicates
the presence of star subgraphs of degree at least $k$ in the correlation network of threshold value $\rho$.
Thus $V_k(\mathbb{X})$ is an extreme value statistic that is only sensitive to
the topology of the correlation network through the distribution of star subgraphs.

\medskip
\begin{theorem}\label{thm:CorrScr}
Let $\mathbf{\Sigma}$, the population dispersion matrix of the rows of $\mathbb{X}$, be row sparse of degree $s=o(p)$. Fix $k$ positive integer, let $p \to \infty$ and $\rho = \rho_p \to 1$ such that
$p^{1/k} (p-1)(1-\rho^2)^{(n-2)/2} \to e_{n,k} \in (0, \infty)$. 
Then:
\begin{enumerate}
\item as $p\rightarrow \infty$:
$$\Prob(V_k(\mathbb{X}) \geq \rho) \to 1-\exp(-\Lambda J_\mathbf{{X}}/\phi(k)),$$
where
$J_\mathbf{X}$ is is a fixed scalar depending on the distribution of the sample correlation matrix $\mathbf R$ (defined in 
\cite[Equation (33)]{hero-bala-IT-2012}) 
and
$$\Lambda = \lim_{p\to \infty, \rho \to 1} \Lambda(\rho) = ((e_{n,k} a_n ) / (n-2) )^{k}/k!,$$
with
$$\Lambda(\rho) = p {p-1 \choose k} P_0(\rho)^k,$$
$$P_0(\rho) = a_n \int_{\rho}^1 (1-u^2)^{\frac{n-4}{2}} du,$$
$$a_n = 2B((n-2)/2, 1/2) \mbox{ with } B(l,m) \mbox{ the beta function}, $$
$\phi(k)=2$ if $k = 1$, and $\phi(k)=1$ otherwise.
\item If $\mathbf{\Sigma}$ is block sparse of degree $s$,
then
$$J_\mathbf{{X}}=1 + O((s/p)^{k+1}),$$
and if $\mathbf{\Sigma}$ is diagonal then $J_\mathbf{{X}}=1$.
\end{enumerate}
\end{theorem}
\begin{IEEEproof}
The result follows from \eqref{eq:V_eq_N} and from Proposition 2 in \cite{hero-bala-IT-2012}:
$$\Prob(N_{k, \rho} > 0) \to 1-\exp(-\Lambda J_\mathbf{{X}}/\phi(k)),$$
under the same asymptotic limit of $p$ and $\rho$ specified in the theorem statement.
\end{IEEEproof}
\medskip

In Section~\ref{sec:InterpretJ} below we provide some insights into the nature of the parameter $J_\mathbf{{X}}$ appearing
in Theorem~\ref{thm:CorrScr} above. We comment on the consequences of this theorem.
Furthermore, it can be shown that the Theorem~\ref{thm:CorrScr} holds if $p$ and $n$ both go to infinity as long as $e_{n,k}$ is finite and nonzero. This is guaranteed if $n$ increases as $\log p$. 

Using \eqref{eq:V_eq_N} and Theorem~\ref{thm:CorrScr},
the  large $p$  distribution of $V_k$ defined in \eqref{eq:Sum_Stat_V} can be approximated by
\begin{equation}\label{eq:CDF_V}
\Prob(V_k(\mathbb{X}) \leq \rho) = \exp(-\Lambda(\rho) J_\mathbf{{X}}/\phi(k)), \; \rho \in [0,1],
\end{equation}
where $\Lambda(\rho)$ is as defined in Theorem~\ref{thm:CorrScr}.
Although the limits in Theorem~\ref{thm:CorrScr} are guaranteed to hold for large values of $\rho$, numerical experiments \cite{hero-bala-IT-2012}
have shown that the approximation \eqref{eq:CDF_V} remains accurate for smaller values of $\rho$ as long as $n$ is small and $p\gg n$.

The distribution \eqref{eq:CDF_V} is differentiable everywhere except at $\rho=0$
since $P(V_k(\mathbb{X})=0)>0$. For $\rho>0$ and large $p$,  $V_k$ has density
\begin{equation}\label{eq:PDF_V_1}
f_V(\rho) = - \frac{\Lambda'(\rho) }{\phi(k)}J_\mathbf{{X}} \exp\left(-\frac{\Lambda(\rho)}{\phi(k)} J_\mathbf{{X}}\right), \; \rho \in (0,1].
\end{equation}
Note that $f_V$ in \eqref{eq:PDF_V_1} is the density of the Lebesgue continuous component of
the distribution \eqref{eq:CDF_V} and that it integrates over $\rho\in (0,1]$ to $1-O(e^{-p^2})$.

The density $f_V$ is a member of a one-parameter exponential family with $J_\mathbf{X}$ as the unknown parameter.
This follows from the following relations. First
\begin{equation}\label{eq:Lambdarho}
\begin{split}
\Lambda(\rho) &=  p {p-1 \choose k} \left(a_n \int_{\rho}^1 (1-u^2)^{\frac{n-4}{2}} du\right)^k \\
&= C \; T(\rho)^k,
\end{split}
\end{equation}
where
\begin{equation}\label{eq:C_pndelta}
C=C_{p,n,k}=p {p-1 \choose k} a^k_n
\end{equation}
does not depend on $\rho$, and
\begin{equation}\label{eq:Trho}
T(\rho) = \int_{\rho}^1 (1-u^2)^{\frac{n-4}{2}} du
\end{equation}
 is the incomplete beta function.
Using \eqref{eq:Lambdarho} and noting that $T(\rho)' = -(1-\rho^2)^{\frac{n-4}{2}}$, $f_V(\rho) = f_V(\rho; J_\mathbf{{X}})$
is a member of the exponential family with parameter $J_\mathbf{{X}}>0$:
\begin{equation}\label{eq:PDF_V}
\begin{split}
f_V &(\rho; J_\mathbf{{X}})  \\
=& \frac{C k}{\phi(k)} T(\rho)^{k-1} (1-\rho^2)^{\frac{n-4}{2}} J_\mathbf{{X}} \exp\left(-\frac{C  T(\rho)^k}{\phi(k)} J_\mathbf{{X}}\right).
\end{split}
\end{equation}

The vertex degree parameter $k$ in \eqref{eq:PDF_V} is a fixed design parameter
that can be selected to maximize change detection performance according to \eqref{prob:Pollak}.
In the sequel, we focus on the case $k=1$ for simplicity of presentation. 
All of the analysis below continues to hold when $f_V$ is replaced by the general form \eqref{eq:PDF_V}. 
For $k=1$, the statistic $V_k$ reduces to the maximal magnitude correlation
\begin{equation}\label{eq:Sum_Stat_V_deltaeq1}
V(\mathbb{X})= \max_{i \neq j} |\mathbf{R}_{ij}|,
\end{equation}
and the density in \eqref{eq:PDF_V} reduces to
\begin{equation}\label{eq:PDF_V_delta1}
f_V(\rho; J) = \frac{C}{2} (1-\rho^2)^{\frac{n-4}{2}} J \exp\left(-\frac{C}{2} J \;T(\rho) \right), \; \rho \in (0,1],
\end{equation}
where we have suppressed subscript $\mathbf{X}$ in the exponential family parameter $J$. 

\begin{figure}[htb]
\center \vspace{-0.3cm}
\includegraphics[width=8cm, height=5cm]{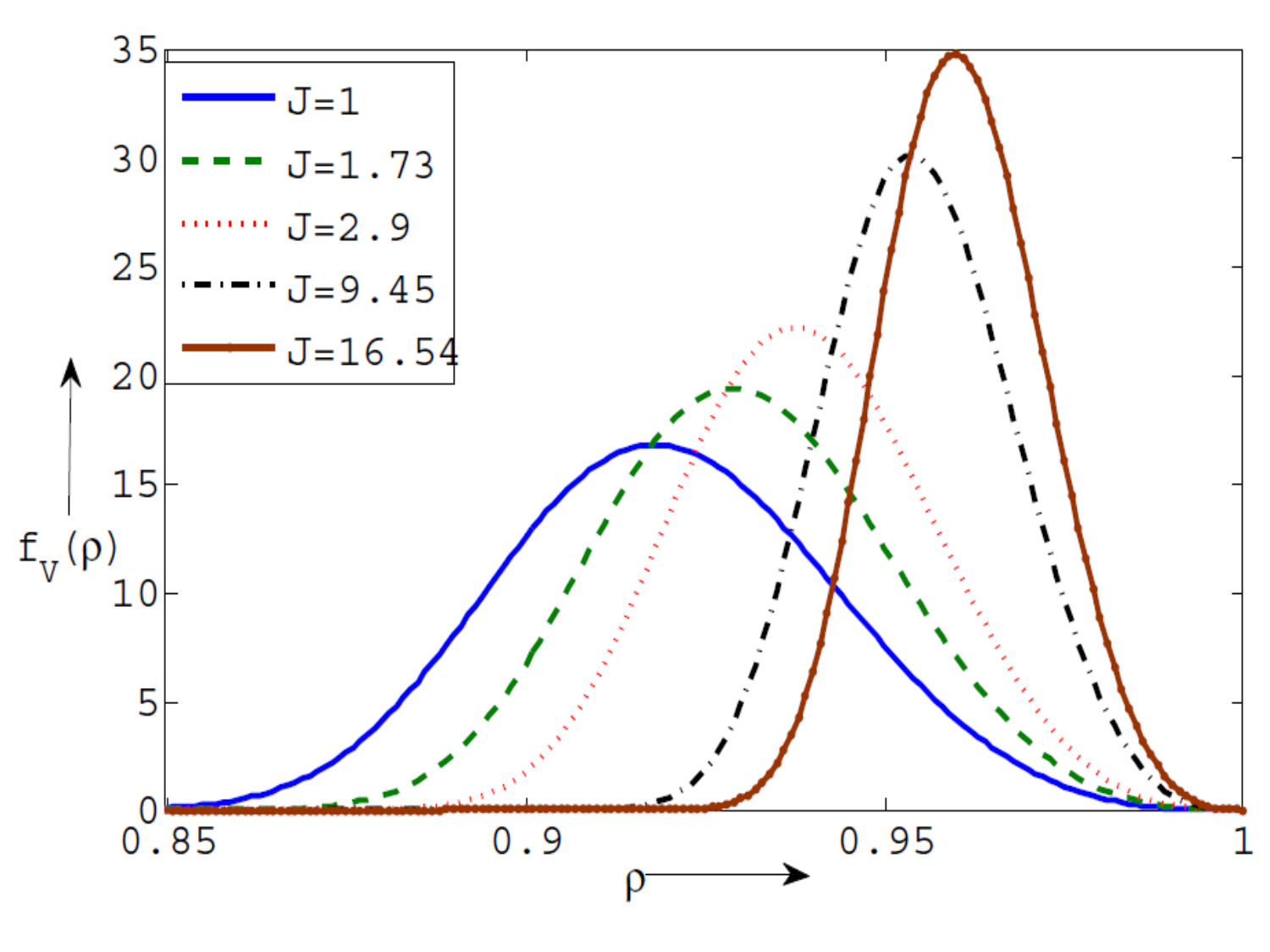}\vspace{-0.5cm}
\caption{Plot of density $f_V$ in \eqref{eq:PDF_V} for various values of the parameter $J$ for $n=10, p=100$.
 This is the density of the summary statistic used to detect a change in covariance of the random matrix sequence $\mathbb X$.}
\label{fig:fV_density}\vspace{-0.3cm}
\end{figure}

In Fig.~\ref{fig:fV_density} is plotted the density $f_V$ for various values of $J$ for
$n=10$, and $p=100$. We note that for the chosen values of $n$ and $p$, the density is concentrated close to $\rho=1$.

\subsection{Interpretation of the parameter $J=J_\mathbf{{X}}$}
\label{sec:InterpretJ}
The asymptotic approximation to the probability $\Prob(N_{k, \rho} > 0)$,
used in Theorem~\ref{thm:CorrScr}, is obtained in \cite{hero-bala-IT-2012} by relating
$N_{k, \rho}$ to a Poisson random variable in the purely high dimensional limit as $p \rightarrow \infty$ and $n $ fixed.
The first step is to recall the Z-score representation of the sample correlation $\bR$:
$$\bR=\mathbb Z^T\mathbb Z, \hspace{0.2in} \mathbb Z=[\bZ_1, \ldots,\bZ_p] $$
$$\bZ_i= \frac{\bX_i-\hat{\bmu}_i \mathbf 1}{\sqrt{\hat{\sigma}}_{ii}\sqrt{n-1}} , \;\; i=1,\ldots, p$$
These Z-scores lie in a $n-2$ dimensional subspace
$$\mathbf1^T \bZ_i  =0 \; \mbox{and} \;  \|\bZ_i\| =1.$$
Due to the fact that
$$
\bZ_i^T\bZ_j=\bR_{ij}, \;\; \mbox{and} \;\; \|\bZ_i - \bZ_j\| = \sqrt{2(1-\bR_{ij})}, $$
the correlation between the columns of the data matrix is directly related to the
Euclidean distance between their corresponding Z-scores. The parameter
$J=J_{\mathbf{X}}$ is a limiting value of an average of the joint density of the Z-scores.
It is a Z-score uniformity measure: $J=1$ implies the scores are uniformly distributed
on the $n-2$ dimensional sphere, $J>1$ if the scores are homophilic in nature, and $J<1$ if they are homophobic.
For more details see Section-II in \cite{hero-bala-jasta-2011}.

\section{QCD for large scale random matrices}\label{sec:QCD_on_V}
In this section we apply the asymptotic results of Theorem~\ref{thm:CorrScr} to quickest
change detection of the maximal kNN coherence in the data matrix sequence $\{\mathbb{X}(m)\}$.
We assume that both the pre- and post-change dispersion matrices,
$\mathbf{\Sigma}_0$ and $\mathbf{\Sigma}_1$, are row sparse with degree $s=o(p)$,
and map the data matrix sequence $\{\mathbb{X}(m)\}$ to the sequence of summary statistics $\{V_k(\mathbb{X}(m))\}_{m \geq 1}$,
with $k=1$.
For simplicity we refer to this sequence as $\{V(m)\}$.
In the asymptotic regime considered in Theorem~\ref{thm:CorrScr}, the random variables $\{V(m)\}$
each have an approximate asymptotic density $f_V(\cdot; J)$ of form \eqref{eq:PDF_V_delta1}. 
Let $J_0$ and $J_1$ be the value of parameter $J$ before and after change point $\gamma$, 
respectively.
The QCD problem is to detect and localize the change in distribution of $V(m)$:
\begin{equation}\label{eq:QCD_in_V}
\begin{split}
V(m)     &\sim f_V(\cdot; J_0), \; m < \gamma \\
 V(m)    &\sim f_V(\cdot; J_1), \; m \geq \gamma.
\end{split}
\end{equation}
Below we consider two cases of row-sparse pre-change dispersion matrices:
\begin{enumerate}
\item $\mathbf{\Sigma}_0$ is diagonal, and
\item $\mathbf{\Sigma}_0$ is not diagonal but is block-sparse.
\end{enumerate}
Note that $\mathbf{\Sigma}_1$ is only assumed to be row-sparse.

If $\mathbf{\Sigma}_0$ is diagonal then, from Theorem~\ref{thm:CorrScr}, $J_0=1$, and the QCD problem
in \eqref{eq:QCD_in_V} reduces to detecting a change in parameter $J$ of the exponential family density  \eqref{eq:PDF_V_delta1},
with \textit{known} pre-change parameter value.
The change in this case can be efficiently detected using Lorden's GLR test \cite{lord-amstat-1971}
(also see Section~\ref{sec:knownJ0} below), and the
test can be designed using the performance analysis provided in \cite{lord-amstat-1971}.

In the case of non-diagonal dispersion matrix $\mathbf{\Sigma}_0$, the QCD problem in \eqref{eq:QCD_in_V}
has an \textit{unknown} pre-change parameter $J_0$. There are no known efficient solutions
to the QCD problem in this case\footnote{The difficulty in this setting is that due to our minimax formulation, depending on the algorithm used, the worst case delay for the chosen algorithm will occur when we do not have enough time 
to learn the pre-change parameter.}. 

However, we recall that
if the dispersion matrix $\mathbf{\Sigma}_0$ is only block-sparse with degree $s \ll p$,
then by assertion 2 of Theorem~\ref{thm:CorrScr}, $J_0$ is close to $1$, i.e., $|J_0-1|$ is small.
Motivated by this fact we use Lorden's test in this case as well, with $J_0$ set equal to $1$,
and characterize the range of pre-change parameters close to $1$ for which the change can be detected efficiently.
Specifically, in Section~\ref{sec:unknownJ0} below, we provide delay and false alarm analysis of Lorden's test
when there is a mis-specification in the pre- and post-change distribution. 
As discussed in the introduction, such an analysis for SPRT and CUSUM is 
standard (though nontrivial), and is available in the literature \cite{tart-niki-bass-2014}. But, the one involving GLR based CUSUM, i.e., Lorden's test 
requires an extra set of conditions. These extra set of conditions are specified in 
Assumptions~\ref{assum:kappa_gtheta}--\ref{assump:I_theta} below. 

We note again that the performance analysis in Section~\ref{sec:unknownJ0} is provided for an arbitrary
one-parameter exponential family, and not just for the family in \eqref{eq:PDF_V_delta1}.

\subsection{QCD with Diagonal Pre-Change Dispersion Matrix $\mathbf{\Sigma}_0$}
\label{sec:knownJ0}
If the pre-change dispersion matrix is diagonal, then from Theorem~\ref{thm:CorrScr} $J_0=1$,
and the QCD problem \eqref{eq:QCD_in_V} reduces to the
parametric QCD problem with unknown post-change parameter $J$:
\begin{equation}\label{eq:QCDProb_1toother}
\begin{split}
V(m) &\sim f_V(\cdot; 1), \quad\hspace{1.25cm} m < \gamma \\
    &\sim f_V(\cdot; J), \; \quad J \neq 1, \; m \geq \gamma.
\end{split}
\end{equation}

Consider the following QCD test, defined by the stopping time $\taug$\footnote{The subscript G in $\taug$ is used to denote a GLR test. This is not to be confused with the use of 
density function $g$ in the misspecification analysis to follow.}
\begin{equation}\label{eq:taug}
\begin{split}
\taug =\inf_{m \geq 1}  \left\{\max_{1 \leq \ell \leq m} \sup_{J \in \mathcal{J}: |J-1| \geq \epsilon} \sum_{i=\ell}^m \log \frac{f_V(V(i);J)}{f_V(V(i); 1)} > A\right\},
\end{split}
\end{equation}
where $A$ and $\epsilon>0$ are user-defined parameters, and $\mathcal{J}$ is a user-defined set of post-change parameters. 
The parameter $A$ is a threshold used to control the false alarm rate.
The parameter $\epsilon$ represents the minimum magnitude of change, relative to $J=1$, that the user wishes to detect.
In the following, we chose either $\mathcal{J} = (-\infty, \infty)$ or $\mathcal{J} = [0, \infty)$. 

The stopping rule $\taug$ was shown to be asymptotically optimal in \cite{lord-amstat-1971} for a related QCD problem
when 1) the marginal density $f_V(v;\cdot)$ of the observation sequence $\{V(m)\}$
is of known form that is a member of a one-parameter exponential family, and 2) when the parameter $J_0$  of the pre-change density is known.
Both of these properties are satisfied for the summary statistic $V=V(\mathbb X)$ for
the QCD model in \eqref{eq:QCDProb_1toother}, since $J_0=1$.
Due to the results in \cite{lai-ieeetit-1998}, the stopping rule $\taug$ is asymptotically
optimal for the minimax change detection problem \eqref{prob:Pollak} as well.

The following theorem establishes strong asymptotic optimality of this test with $\{V(m)\}$ as
the observation sequence. It also provides delay and false alarm estimates with which the test can be designed.
\medskip
\begin{theorem}[\cite{lord-amstat-1971}, \cite{lai-ieeetit-1998}]\label{thm:Lorden}
Fix any $\epsilon > 0$ and $\mathcal{J} = (-\infty, \infty)$. 
\begin{enumerate}
\item For the stopping rule $\taug$, the supremum in \eqref{prob:Pollak} is
achieved at $\gamma=1$, i.e.,
\[
\sup_{\gamma\geq 1} \; \Expect_\gamma[\taug-\gamma| \taug \geq \gamma] = \Expect_1[\taug-1].
\]
\item For $C_\epsilon > 0$ depending on $\epsilon$, setting $A=\log [C_\epsilon \beta \log \beta ]$ ensures that as $\beta \to \infty$,
\[
\Expect_\infty[\taug] \geq \beta (1+o(1)),
\]
and for each possible true post-change parameter $J$, with $|J-1|\geq \epsilon$,
\begin{equation}
\begin{split}
\Expect_1[\taug] &= \frac{\log \beta}{I(J)} (1+o(1)) \\
&= \inf_{\tau: \Expect_\infty[\tau] \geq \beta}\sup_{\gamma\geq 1} \; \Expect_\gamma[\tau-\gamma| \tau \geq \gamma] (1+o(1)),
\end{split}
\end{equation}
where $I(J)$ is the Kullback-Leibler divergence between the densities $f_V(\cdot; J)$ and $f_V(\cdot; 1)$.
\end{enumerate}
\end{theorem}
\medskip
Theorem~\ref{thm:Lorden} implies that the stopping rule $\taug$ is uniformly asymptotically optimal for each post-change parameter $J$,
as long as $|J-1| \geq \epsilon$. For convenience of implementation one can also use the window limited variation
of $\taug$ as suggested in \cite{lai-ieeetit-1998}. The constant $C_\epsilon$ in the theorem can be explicitly obtained from Theorem~\ref{thm:Misspec_FAR} below; see equations \eqref{eq:LBagreementeqn} and \eqref{eq:Cepsilonvalue} after the theorem. Also note that, the theorem above is true for any $\mathcal{J} \subset (-\infty, \infty)$ since among all choices of parameter set $\mathcal{J}$, $\mathcal{J} = (-\infty, \infty)$ leads to the worst case trade-off between the delay and the rate of false alarms. 

\subsection{QCD with Block-Sparse Pre-Change Dispersion Matrix $\mathbf{\Sigma}_0$ }
\label{sec:unknownJ0}

As discussed above, if the pre-change dispersion matrix $\mathbf{\Sigma}_0$ is not diagonal,
then the pre-change parameter $J_0\neq 1$. If $\mathbf{\Sigma}_0$ is block-sparse with degree $s$,
then by part 2 of Theorem~\ref{thm:CorrScr}, $|J_0-1|$ is small.
This motivates the use of Lorden's test as in \eqref{eq:taug} for QCD. However, Theorem~\ref{thm:Lorden} no longer applies 
since $f_V(\cdot; J_0)$ with $J_0=1$ is a mis-specification of the true pre-change distribution.
In this section we extend Theorem~\ref{thm:Lorden} to cover stopping rules specified by  Lorden's GLR test for $\{f_V(\cdot; J)\}$ vs $f_V(\cdot; 1)$ under this kind of mismatch. The theorem proven below is in fact applicable to a broader class of scalar parameter exponential families,
not just to the $\{f_V(\cdot; J)\}$ family \eqref{eq:PDF_V_delta1} considered in this paper.

Consider the following general setting. For scalar parameter $\theta$ let $\{f_\theta\}$ be a parametric exponential family of distributions with respect to a $\sigma$-finite measure $\mu$
\begin{equation}
\label{eq:t_expo}
f_{\theta}(y) = e^{\theta T(y) - b(\theta)} h(y), \quad \theta \in \Theta,
\end{equation}
where $\Theta$ is a specified interval on the real line and $b(\theta)$ is differentiable everywhere on $\Theta$.

The QCD test $\taug$ for
detecting a change from $f_{\theta_0}$ to $f_\theta$, under the constraints $\theta \in \Theta$ and $|\theta-\theta_0| \geq \epsilon$ is given by
\begin{equation}
\label{t_GLR}
\taug = \inf\left\{m\geq 1: \max_{1\leq k \leq m} \; \sup_{\theta: |\theta-\theta_0| \geq \epsilon} \; \sum_{i=k}^m \log \frac{f_\theta(Y_i)}{f_{\theta_0}(Y_i)} > A\right\},
\end{equation}
where $\{Y_i\}$ is an i.i.d. observation sequence. The $\taug$ given in \eqref{eq:taug} is a special case of the $\taug$ given in \eqref{t_GLR} 
with $f_\theta(\cdot)$ replaced by $f_V(\cdot; J)$. 
Below we provide performance bounds for the mean time to false alarm and
the average detection delay when the samples are drawn from a density $g$. Specifically,
we provide bounds on  $\Expect_g[\taug]$, where $\Expect_g$ denotes expectation
with respect to the probability measure under which all the samples $\{Y_i\}$
have density $g$. When the density $g$ is close to $f_{\theta_0}$, the expression $\Expect_g[\taug]$
can be interpreted as an estimate of the mean time to false alarm. When the density $g$ is close to $f_{\theta}$,
for some $\theta$ with $|\theta-\theta_0| \geq \epsilon$, then the expression $\Expect_g[\taug]$
can be interpreted as an estimate of the average time to change detection.

\medskip
\subsection{Mean Time to False Alarm}
We first provide a lower bound on $\mathsf{E}_g[\taug]$ when $g$ is not necessarily equal to $f_{\theta_0}$,
but is close to $f_{\theta_0}$ in a particular sense.
The closeness is characterized through the following three assumptions.

  \smallskip
\begin{assumption}\label{assum:kappa_gtheta}
There exists a positive constant $\kappa_{\theta, g}$ such that for every $\theta \in \Theta$
with $|\theta-\theta_0| \geq \epsilon$
\begin{equation}\label{eq:kappaCond}
\int \left(\frac{f_\theta(y)}{f_{\theta_0}(y)}\right)^{\kappa_{\theta, g}} g(y) \; d\mu(y)= 1.
\end{equation}
Furthermore, there exists $\kappa_g$ such that
\begin{equation}\label{eq:kappastar}
0 < \kappa_g \leq \inf_{\theta \in \Theta: |\theta-\theta_0| \geq \epsilon} (\kappa_{\theta, g} ).
\end{equation}
\end{assumption}
\smallskip

The condition in \eqref{eq:kappaCond} is the classical condition needed to analyze one-sided tests under mis-specification \cite{wood-nonlin-ren-th-book-1982}.
The condition in \eqref{eq:kappastar} is an additional condition that will be needed for analysis of the GLR test  defining the stopping time \eqref{t_GLR}.

Let $\mathcal{G}$ be a family of densities on the real line, for example,
$\mathcal{G} \subset \{f_\theta\}$.
\smallskip
\begin{assumption}\label{assump:kappa_g}
There exists a positive constant $\kappa^*$ such that
\begin{equation}
0 < \kappa^* \leq \kappa_g, \quad \forall g \in \mathcal{G}.
\end{equation}
\end{assumption}
\smallskip

\smallskip
\begin{assumption}\label{assump:I_theta}
The KL-divergence $\text{KL}(f_\theta \| f_{\theta_0})=I(\theta)$ between $f_\theta$ and $f_{\theta_0}$ increases with $|\theta-\theta_0|$.
\end{assumption}
\medskip

\begin{theorem}\label{thm:Misspec_FAR}
\begin{enumerate}
\item
If Assumption~\ref{assum:kappa_gtheta} and Assumption~\ref{assump:I_theta} are satisfied then
$$\mathsf{E}_g[\taug] \geq \frac{e^{\kappa_{g} A}}{2 \left( \frac{A}{\min \{ I(\theta_0 + \epsilon),  I(\theta_0 - \epsilon)\}}+1  \right)}.
$$
\item Furthermore, if Assumption~\ref{assump:kappa_g} is also satisfied then for every $g \in \mathcal{G}$
$$\mathsf{E}_g[\taug] \geq \frac{e^{\kappa^* A}}{2 \left( \frac{A}{\min \{ I(\theta_0 + \epsilon),  I(\theta_0 - \epsilon)\}}+1  \right)}.
$$
\end{enumerate}
\end{theorem}
\begin{IEEEproof}
See appendix.
\end{IEEEproof}

\medskip
We note that the lower bound in the second part of the above theorem is not a function of the density $g$.
Also, if $g=f_{\theta_0}$ then $\kappa^* = \kappa_g = \kappa_{\theta, g}=1$, $\forall \theta$,
then the lower bounds agree with the bounds presented in Theorem~\ref{thm:Lorden}. Specifically, by setting $\kappa^* = \kappa_g =1$ and $A= \log [C_\epsilon \beta \log \beta ]$ we get as $\beta \to \infty$
\begin{equation}\label{eq:LBagreementeqn}
\begin{split}
\mathsf{E}_g[\taug] &\geq \frac{e^{A}}{2 \left( \frac{A}{\min \{ I(\theta_0 + \epsilon),  I(\theta_0 - \epsilon)\}}+1  \right)}\\
                                &= \frac{C_\epsilon \beta \log \beta}{2 \left( \frac{\log [C_\epsilon \beta \log \beta ]}{\min \{ I(\theta_0 + \epsilon),  I(\theta_0 - \epsilon)\}}+1  \right)}\\
                                &=\frac{C_\epsilon \beta \log \beta}{2 \left( \frac{\log \beta }{\min \{ I(\theta_0 + \epsilon),  I(\theta_0 - \epsilon)\}}\right)}(1+o(1))\\
                &=\frac{C_\epsilon \beta}{\frac{2}{\min \{ I(\theta_0 + \epsilon),  I(\theta_0 - \epsilon)\}}}(1+o(1)).                
\end{split}
\end{equation}
Thus, the right choice of $C_\epsilon$ is 
\begin{equation}\label{eq:Cepsilonvalue}
C_\epsilon = \frac{2}{\min \{ I(\theta_0 + \epsilon),  I(\theta_0 - \epsilon)\}}.
\end{equation}

\medskip
\subsection{Average Detection Delay} 
We next obtain an upper bound on $\Expect_g[\taug]$ when $g$ is close
to one of the members of the post-change set of densities  $\{f_\theta: \theta \in \Theta; |\theta-\theta_0|>\epsilon\}$. The closeness here is characterized by the following assumption.
\begin{assumption}\label{assump:delay}
 $\exists \theta_{g}$, s.t. $|\theta_g-\theta_0| \geq \epsilon$ \quad
 $$\int \; \log [f_{\theta_{g}}(y)/f_{\theta_0}(y)] \; g(y) \; d \mu(y) > 0.$$
\end{assumption}
\medskip

\begin{theorem}\label{thm:Misspec_Delay}
If Assumption~\ref{assump:delay} is satisfied then
$$
\Expect_g[\taug] \leq \frac{A}{\int \log [f_{\theta_{g}}(y)/f_{\theta_0}(y)] g(y) d\mu(y) }(1+o(1)) \mbox{ as } A \to \infty.
$$
\end{theorem}
\begin{IEEEproof}
See appendix.
\end{IEEEproof}
\medskip
We note that if $g=f_\theta$, for $|\theta-\theta_0| \geq \epsilon$, then the upper bound in Theorem~\ref{thm:Misspec_Delay} is
the mean detection delay of the GLR test as obtained in Theorem~\ref{thm:Lorden}.

\vspace{0.2cm}
\subsection{Discussion on the Assumptions}
For the Theorem~\ref{thm:Misspec_FAR} and Theorem~\ref{thm:Misspec_Delay} to be used to bound the mean time to false alarm and the mean detection delay, the distribution family $f_\theta$ must satisfy Assumptions~\ref{assum:kappa_gtheta}--\ref{assump:delay} stated above. In some cases these conditions can be analytically verified, e.g., the case of detecting a shift of mean in the Gaussian distribution as shown below.

Consider the Gaussian density parameterized by its mean $\theta$
\begin{equation}\label{eq:GaussDenAssumptions}
f_\theta(y) = \frac{1}{\sqrt{2 \pi}} e^{-\frac{(y-\theta)^2}{2}}, \theta \in (-\infty, \infty), \; y \in (-\infty, \infty)
\end{equation}
where the objective is to detect a change in mean from a level $\theta_0$ to a level $\theta$, with $|\theta - \theta_0| \geq \epsilon$. 
Let the samples $\{Y_m\}$ have density $g=f_{\tilde{\theta}_0}$. 
In this case the integral expression \eqref{eq:kappaCond} has a closed form expression and the following establish that Assumptions~\ref{assum:kappa_gtheta}--\ref{assump:I_theta} are satisfied. The following is proven in the Appendix.

\begin{lemma}\label{lem:Gaussian}
Let $f_\theta$ be of the form \eqref{eq:GaussDenAssumptions} and fix $\theta_0$, $\tilde{\theta}_0$. Then 
\begin{enumerate}
\item The KL divergence KL$(f_\theta \| f_{\theta_0})$ is strictly increasing in $|\theta-\theta_0|$. Furthermore, there exists a $\kappa_{\theta, g}>0$ given by
\begin{equation}\label{eq:gauss_kap_g_th}
\kappa_{\theta, g} = 1 + \frac{2(\theta_0 - \tilde{\theta}_0)}{\theta-\theta_0}
\end{equation}
satisfying \eqref{eq:kappaCond} provided $|\tilde{\theta}_0-\theta_0|< \epsilon/2$.
\item There exists a $\kappa_g >0$ given by
\begin{equation}\label{eq:gauss_kap_g}
\kappa_g = \min\{\kappa_{g,\theta_0+\epsilon}, \; \kappa_{g,\theta_0-\epsilon}\} = 1 - \frac{2(|\tilde{\theta}_0-\theta_0|)}{\epsilon}
\end{equation}
that satisfies \eqref{eq:kappastar}.
\item Let
\[
\mathcal{G} = \{f_\theta: |\theta-\theta_0| \leq \epsilon/3\}.
\]
Then
\[
\kappa^* = \min\{\kappa_{f_{\theta_0+\epsilon/3},\theta_0+\epsilon}, \; \kappa_{f_{\theta_0-\epsilon/3},\theta_0-\epsilon}\}
\]
satisfies Assumption~\ref{assump:kappa_g}. That is, $\kappa^*$ is the minimum of $k_g$ with $g=f_{\theta_0+\frac{\epsilon}{3}}$
and $k_g$ with $g=f_{\theta_0-\frac{\epsilon}{3}}$.
\end{enumerate}
\end{lemma}
\medskip

\begin{figure}[htb]
\center \vspace{-0.3cm}
\includegraphics[width=8cm, height=5cm]{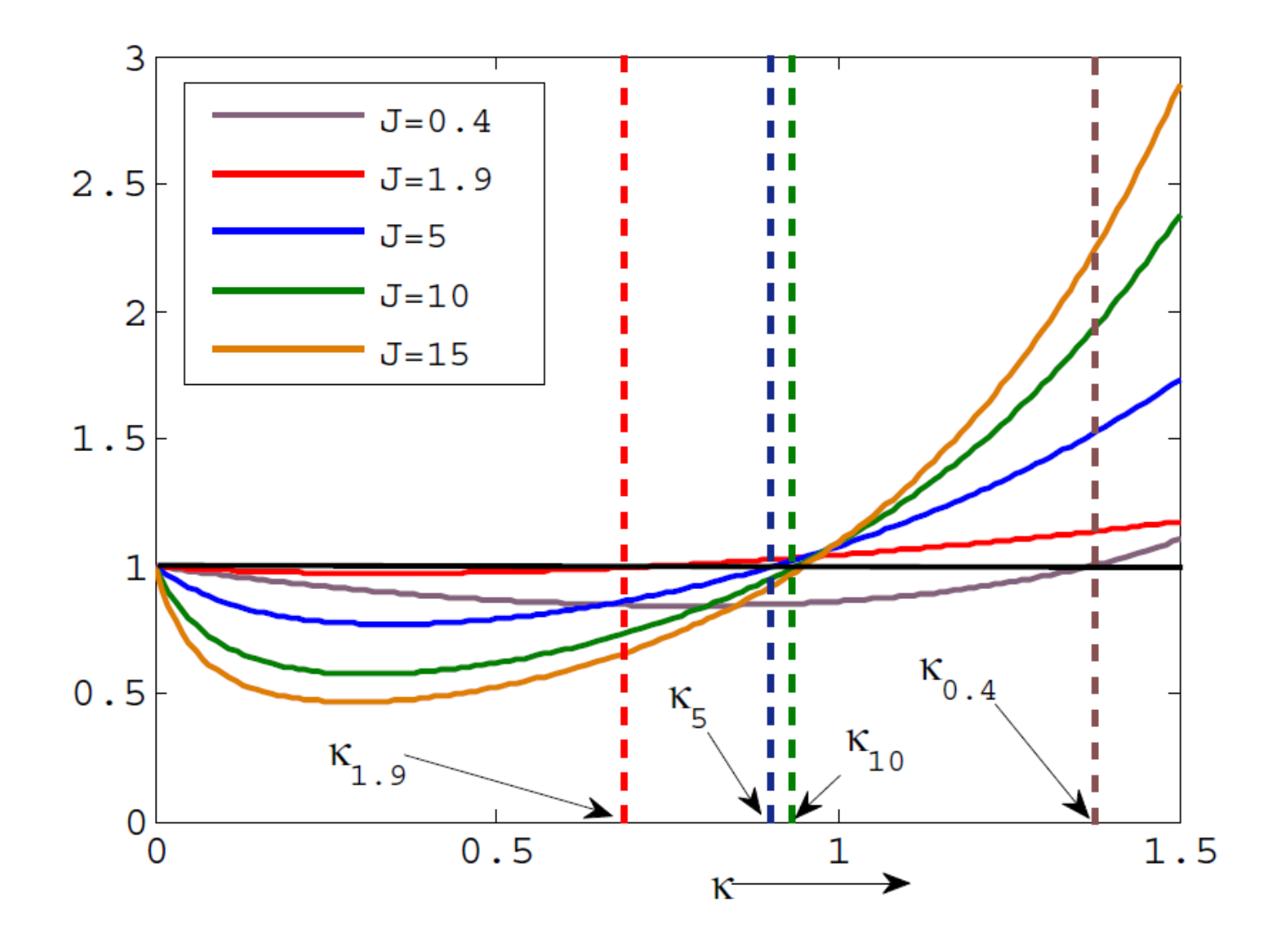}
\caption{Plot of the integral $\int (f_V(v; J)/f_V(v; 1))^{\kappa} f_V(v; J_0) \; dv= 1$ as a function of $\kappa$,
for various values of $J$: $J=0.4, 1.9, 5,10,15$. The dashed lines show the point of intersection of the curve corresponding
to a particular value of $J$ with the straight line at height $1$. Note that the value of $\kappa_g$ of $\kappa$ at which the curves take value $1$ increases with the parameter value $J$ when $J>1$, and $\kappa_g \leq 1$ when $J \geq 1$ and $\kappa_g \geq 1$ for $J < 1$. Thus the smallest $\kappa$ is achieved by the parameter $J=1+\epsilon$, which in this case, is $1+\epsilon=1.9$.}
\label{fig:fV_misspec}
\end{figure}

We turn to the problem of detection of a change in parameter $J$ of the distribution $f_V(\rho;J)$ in  \eqref{eq:PDF_V_delta1} 
of the summary statistic $V(m)$. Unlike the Gaussian case, the KL divergence is not in closed form and the assumptions in Theorem~\ref{thm:Misspec_FAR} can only be verified by numerical evaluation. 
Specifically, we verify that if the samples are drawn from $J_0$, i.e., $g=f_{J_0}$,
then, similar to the Gaussian case, the worst case $\kappa_g$ is achieved at the boundary.
Fixing $\epsilon=0.9$ and $J_0=1.1$, we plot the integral in \eqref{eq:kappaCond} for various values of $J$,
$|J-1| \geq \epsilon$, in Fig.~\ref{fig:fV_misspec}. Figure 2 shows several curves, indexed by $J=0.4, 1.9, 5,10, 15$,
that are numerical evaluations of the integral $\int (f_V(v; J)/f_V(v; 1))^{\kappa} f_V(v; J_0) \; dv= 1$ plotted as a function of $\kappa$.
In the figure we can see the points at which the integral equals 1 (identified for example by the labels $\kappa_{1.9}$ and $\kappa_5$),
and the smallest such point correspond to the curve for parameter $J=1+\epsilon=1.9$. 
By inspecting the behavior of these curves for different values of $J_0$ and $\epsilon$ (not shown) one can see 
varying $J_0$ and $\epsilon$
that there is an interval around $1$ for which the corresponding $\kappa_{g,\theta}$ is positive, and the
smallest value is achieved when the post-change parameter  $J$ equals $1+\epsilon$ or $1-\epsilon$. Finally, a plot of the KL divergence 
KL$(f_V(\rho;J)\|f_V(\rho;1))=I(J)$ as a function of $J$ shows that $I(J)$ increases as a function of $|J-1|$, showing that Assumption~\ref{assump:I_theta} is satisfied. 

\section{Numerical Results}\label{sec:Numerical}
Here we apply the stopping rule $\taug$ in \eqref{eq:taug}, with $\mathcal{J} = [0,\infty]$, to the problem of detecting a change
in the distribution when the $\{\mathbb X(m)\}$ are Gaussian distributed random matrices  with i.i.d. rows each having identical covariance matrix $\mathbf \Sigma$. 
The pre-change value $J_0$ of $J$ is equal to one (diagonal covariance matrix) and the post-change value $J_1$ of $J$ is greater than 1. Specifically,  the pre-change covariance is the $p\times p$ diagonal matrix
$\mathbf{\Sigma}_0=\text{diag}(\sigma^2_i)$, where $\sigma_i^2 >0$.
The post-change covariance matrix $\mathbf{\Sigma}_1$ is a row-sparse matrix of degree $s$,
constructed as follows. A $p \times p$ sample from the Wishart distribution is generated and some of
the entries are forced to be zero in such a way that no row has more than $s$ non-zero elements.
Specifically, we retain the top left $s \times s$ block of the matrix. For each row after $s$ rows, 
$s$ consecutive elements after the diagonal term is retained. The rest of the entries are either set to zero, or filled to maintain symmetry. 
Mathematically, for $i$, $s+1 \leq i \leq (p+s)/2$,
all but the diagonal and the $(p+s+1-i)$th element is forced to zero. Each time an entry $(i,j)$ is set to zero,
the entry $(j,i)$ is also set to zero, to maintain symmetry. Finally, a positively scaled diagonal matrix is added to $\mathbf{\Sigma}_1$
to restore its positive definiteness. We set $n=10$, $p=100$, and $s=5$.

This procedure gives a non-diagonal dispersion matrix $\Sigma_1$ that determined the value of $J=J_1$ after the change point. The value can be estimated empirically by drawing repeated samples $\{\mathbb X(m)\}_{m=1}^T$ from the $\mathcal N(0, \mathbf \Sigma_1)$, computing the summary statistics  $\{V(m)\}_{m=1}^T$ and empirically estimating $J$ using the method of \cite{sricharan2011local} or by maximizing $\prod_{m=1}^T f_V(V(m);J)$ using \eqref{eq:PDF_V_delta1}, giving the maximum likelihood estimate MLE.

To implement $\taug$ we have chosen $\epsilon=1.5$,
and we use the maximum likelihood estimator of $J$ which,
as a function of the $m$ samples $(V(1), \cdots, V(m))$ from $f_V(\cdot, J)$, is given by
\begin{equation}\label{eq:MLEst_J}
\hat{J}(V(1), \cdots, V(m)) = \frac{1}{\frac{C}{2}\frac{1}{m} \sum_{i=1}^m T(V(i))}.
\end{equation}
Specifically,
\begin{equation}
\begin{split}
\argmax_{J: J \geq 2.5} \;  \log \sum_{i=\ell}^m & \frac{f_V(V(i);J)}{f_V(V(i); 1)} \\
&= \max\{2.5, \hat{J}(V(\ell), \cdots, V(m))\}.
\end{split}
\end{equation}

\begin{figure}[htb]
\center
\includegraphics[width=8cm, height=6cm]{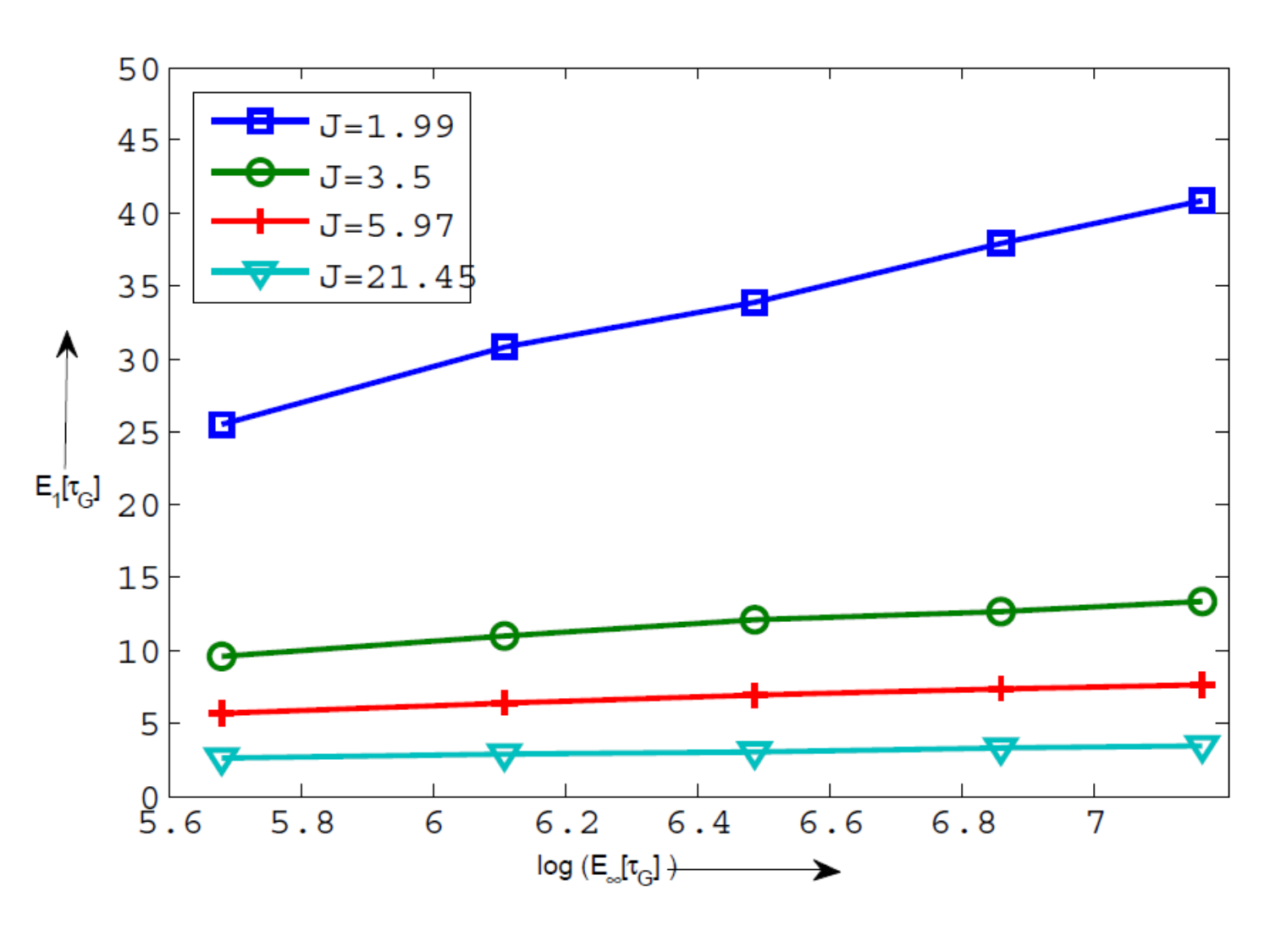}\vspace{-0.2cm}
\caption{The empirical mean time to detect vs  mean time to false alarm (in log scale).
The mean time to detect decreases as the parameter $J$ increases,
 and the relation between $\Expect_1[\taug]$ and $\log (\Expect_\infty[\taug]) $ is linear, as predicted by Theorem 4.1.
 The K-L divergence values for $J=1.99$, $3.5$, $5.97$ and $21.45$ are $0.1906$, $0.5385$, $0.9543$, and $2.1123$, respectively.
 The slopes of the lines are approximately inverse of the K-L divergence values.
}
\label{fig:Delay_MFA}
\end{figure}

In Fig.~\ref{fig:Delay_MFA} we plot the delay to detect ($\Expect_1[\taug]$) vs the log of mean time to false alarm ($\log \Expect_\infty[\taug]$)
for various values of the post-change parameter $J$. The values in the figure are obtained
by choosing different values of the threshold $A$ and estimating the delay by choosing the change point $\gamma=1$
and simulating the test for $500$ sample paths.
The values of the mean time to false alarm are estimated by simulating the test for $1500$ sample paths.
The parameter $J$ for the post-change distribution is estimated
using the maximum likelihood estimator \eqref{eq:MLEst_J}.
As predicted by Theorem~\ref{thm:Lorden}, the delay vs log of false alarm trade-off curve is approximately linear.
For larger values of $J$, the Kullback-Leibler (K-L) divergence between $f_V(\cdot, J)$ and $f_V(\cdot, 1)$ is also larger, resulting
in smaller delays. For the chosen values of the post-change parameters $J=1.99$, $3.5$, $5.97$ and $21.45$,
the corresponding K-L divergence values KL$(f_V(\rho;J)\|f_V(\rho;1))=I(J)$, computed numerically using \eqref{eq:PDF_V_delta1},  are $0.1906$, $0.5385$, $0.9543$, and $2.1123$, respectively.

\begin{figure}[htb]
\center
\includegraphics[width=8cm, height=5cm]{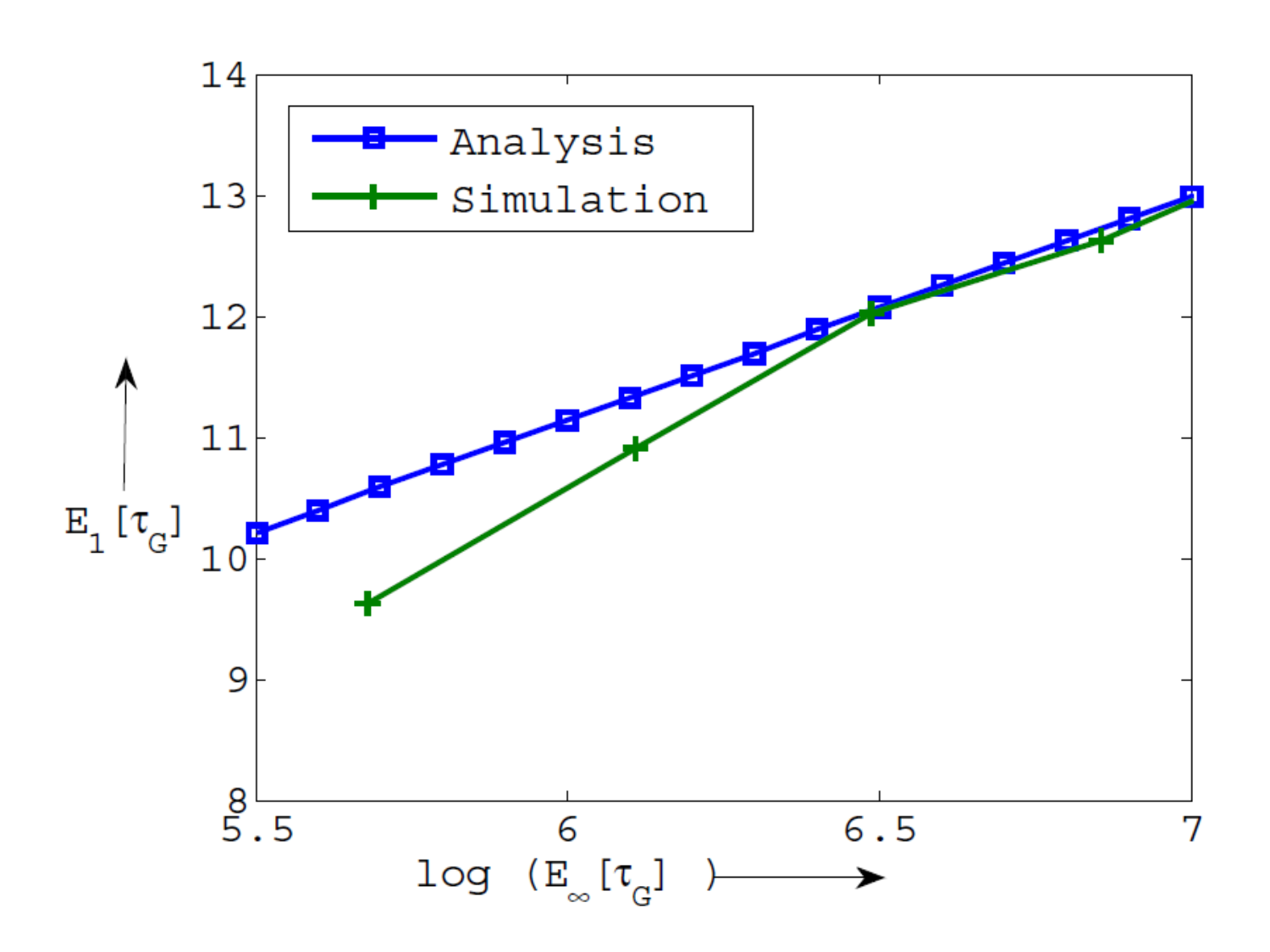}\vspace{-0.3cm}
\caption{Comparison of the delay vs false alarm trade-off curve for $J=3.5$ from Fig.\ref{fig:Delay_MFA} against
the values $\frac{\log \Expect_\infty[\taug]}{I(J)} = \frac{\log \Expect_\infty[\taug]}{0.5385}$ predicted by Theorem~\ref{thm:Lorden}.
As expected, the predictions of the theorem become increasingly accurate as the threshold parameter  $\beta$, and hence the mean time to false alarm, becomes large. }
\label{fig:Delay_MFA_AnaSim}
\end{figure}

In Fig.~\ref{fig:Delay_MFA_AnaSim} we compare the delay vs false alarm trade-off curve
for the post-change parameter $J=3.5$ plotted in Fig.~\ref{fig:Delay_MFA}, against the values $\frac{\log \Expect_\infty[\taug]}{I(J)}$
predicted by Theorem~\ref{thm:Lorden}. Fig.~\ref{fig:Delay_MFA_AnaSim} shows
that the predictions are quite accurate.
We have obtained similar results when the test was simulated for different sparsity degrees $s$.
Thus, the change can be efficiently detected using our proposed methodology.

\begin{figure}[htb]
\center\vspace{-0.3cm}
\includegraphics[width=8cm, height=5cm]{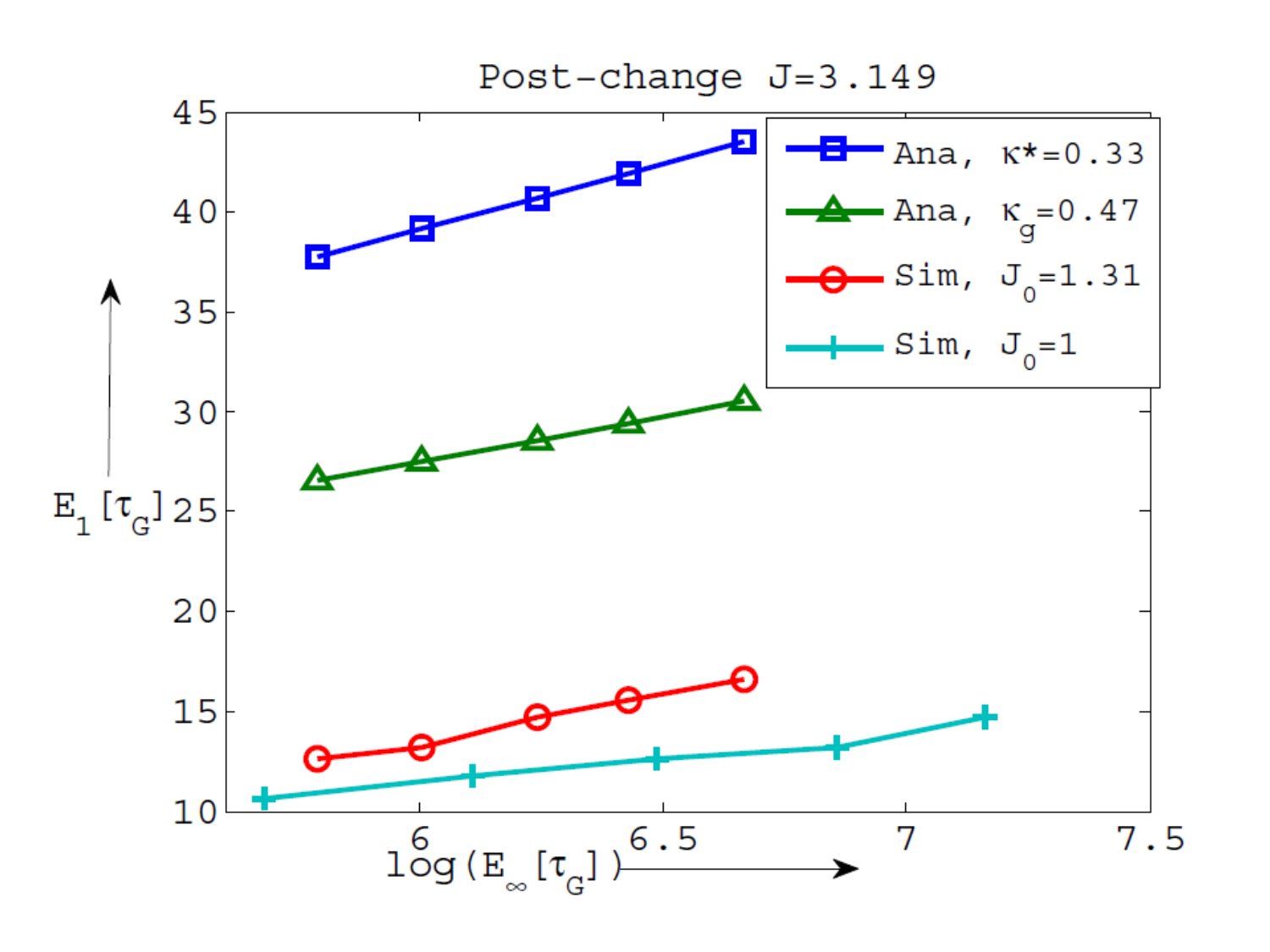}\vspace{-0.3cm}
\caption{Delay-false alarm trade-off curves for Lorden's GLR test of change in distribution of maximal sample coherency $V(m)$ in two misspecification scenarios. Scenario 1 (bottom curve) is the classical regime for Lorden's GLR test, where the pre-change parameter is known ($J_0=1$) but the post-change parameter is unknown (test assumes $J_1>2.5$ while true $J_1$ is $3.149$). In scenario 2 the actual pre-change parameter is $J_0=1.31$, because $\mathbf{\Sigma}_0$ is no longer diagonal (but is block sparse) and the performance of Lorden's GLR test degrades, indicated by the second curve from the bottom. The uncertainty family for the parameter $J_0$ is $\mathcal{G}=\{g=f_{J_0}: |J_0-1| \leq 0.4\}$. 
}
\label{fig:Delay_MFA_AnaSim_Misspec}\vspace{-0.5cm}
\end{figure}

Finally, we evaluate the performance of Lorden's GLR test (15) for detecting changes in maximal coherence when the pre-change or  post-change distributions are possibly misspecified. Fig.~\ref{fig:Delay_MFA_AnaSim_Misspec} shows the delay vs false alarm trade-off curves
for two misspecification scenarios. In both scenarios the experimenter implements Lorden's GLR test assuming that the pre-change parameter is $J=J_0=1$ and the postchange parameter is $J=J_1$ which is unknown except that $J_1\geq 1+\epsilon =1.5$ The first scenario is the standard GLR setting for Lorden's test: the prechange dispersion matrix $\mathbf \Sigma_0$ is diagonal and the post-change matrix $\mathbf \Sigma_1$ is row-sparse giving a true (but unknown) value of $J_1=3.15$. In the second scenario the pre-change dispersion matrix $\mathbf \Sigma_0$ is not diagonal but instead is block-sparse. 
The first curve from the bottom in Fig.~\ref{fig:Delay_MFA_AnaSim_Misspec} is the performance, obtained via simulations, of the Lorden's stopping rule $\taug$ when $\mathbf{\Sigma}_0$ is indeed diagonal and the post-change parameter is $J=3.149$. This curve 
characterizes the performance in the first scenario and will serve as a benchmark  for the second scenario.

In the second scenario, $\mathbf{\Sigma}_0$ is block sparse with block size $5$  and corresponding to parameter value $J=1.31$.
 The performance of Lorden's GLR stopping rule $\taug$ for this case is shown by the second curve from the bottom in Fig.~\ref{fig:Delay_MFA_AnaSim_Misspec}. As expected there is a loss in performance because of misspecification of the pre-change parameter $J_0$. For this plot the threshold $A$ for $\taug$ was chosen
using the knowledge of the pre-change parameter $J_0$.

The remaining (top two) curves in Fig.~\ref{fig:Delay_MFA_AnaSim_Misspec} show
the loss in performance as predicted by Theorem~\ref{thm:Misspec_FAR}. As suggested by the theorem
the asymptotic large $\beta$ delay-false alarm trade-off is given by
\begin{equation}\label{eq:NumRes_kappag}
\Expect_1[\taug] = \frac{\log \beta}{\kappa_g I(J)}(1+ o(1))
\end{equation}
when the pre-change parameter is known, and by
\begin{equation}\label{eq:NumRes_kappastar}
\Expect_1[\taug] = \frac{\log \beta}{\kappa^* I(J)}(1+ o(1))
\end{equation}
when the pre-change parameter is known only within a range of uncertainty.
The top most curve in Fig.~\ref{fig:Delay_MFA_AnaSim_Misspec} is the trade-off curve for the case when
we only know that $|J_0-1| \leq 0.4$.
For this curve, $\kappa^*=0.33$ is obtained by solving $\int (f_J(v)/f_{1}(v))^{\kappa} f_{J_0}(v) \; dv= 1$
for $J=2.5$ and $J_0=1.4$. The second curve from the top is
the trade-off curve when the value of $\kappa$ is obtained by using the knowledge that the pre-change parameter $J_0=1.31$.
The value of $\kappa=\kappa_g$ so obtained is $0.47$.
The top two curves are obtained by dividing the log of mean time to false alarm (the value plotted on the horizontal axis) by the 
lower bounds on mean time to detect, the $o(1)$ approximations on the right hand sides of \eqref{eq:NumRes_kappag} and \eqref{eq:NumRes_kappastar}. 
The fact that these top two curves  are significantly greater than 1 reflects the conservativeness of these lower bounds.

\section{Conclusions and Future Work}
We have introduced and analyzed a method for quickest change detection (QCD) of a step change in the maximal coherence of a sequence of random matrices, under the assumption that the rows of these matrices  are elliptically distributed with row-sparse dispersion matrices before and after the change. In the case that columns of the random matrices are incoherent, i.e., the dispersion matrix is in fact diagonal before the change, the proposed QCD algorithm is first-order asymptotically optimal in the sense of Lorden \cite{lord-amstat-1971} and Pollak \cite{poll-astat-1985}
among all detection algorithms that use the proposed summary statistic, which is the maximal coherence of the sample correlation matrix.
We have also provided mis-specification analysis of the proposed procedure when the pre-change or post-change distributions are unknown.
Future work will include extensions to local summary statistics and experiments
with QCD for sequential dependency testing in high dimension.

\section*{Appendix}
\begin{IEEEproof}[Proof of Theorem~\ref{thm:Misspec_FAR}]
As shown in \cite{lord-amstat-1971}, the key to analysis of $\taug$ is the following
one-sided GLR test
\begin{equation}
\label{t_GLR_SPRT}
N = \inf\left\{m\geq 1: \sup_{\theta: |\theta-\theta_0| \geq \epsilon} \; \sum_{i=1}^m \log \frac{f_\theta(Y_i)}{f_{\theta_0}(Y_i)} > A\right\}.
\end{equation}
Specifically, for any density $g$,
\begin{equation}
\label{eq:onesided_toCuSum}
\Prob_g(N < \infty) \leq \alpha   \quad \implies \quad \Expect_g[\taug] \geq \frac{1}{\alpha},
\end{equation}
where $\Prob_g$ is the probability measure under which all the observations $\{Y_m\}$ have density $g$,
and $\Expect_g$ is the corresponding expectation.
We thus focus on obtaining a bound on $\Prob_g(N < \infty)$.

In reference to this we define the one-sided test between $\theta$ vs $\theta_0$ as
\begin{equation}
\label{t_SPRT}
\nu(f_\theta,f_{\theta_0}) = \inf\left\{m\geq 1: \sum_{i=1}^m \log \frac{f_\theta(Y_i)}{f_{\theta_0}(Y_i)} > A\right\}.
\end{equation}
From Theorem 3.4 in \cite{wood-nonlin-ren-th-book-1982} we know that
if there exists a positive constant $\kappa_{\theta, g} > 0$ such that (see \eqref{eq:kappaCond})
\begin{equation}\label{eq:kappaCond_app}
\int \left(\frac{f_\theta(y)}{f_{\theta_0}(y)}\right)^{\kappa_{\theta, g}} g(y) d\mu(y)= 1,
\end{equation}
then
\begin{equation}\label{eq:UB_SPRT_mis}
\Prob_g(\nu(f_\theta,f_{\theta_0}) < \infty) \leq e^{-\kappa_{\theta, g} A}.
\end{equation}
The basic idea behind \eqref{eq:UB_SPRT_mis} is that if \eqref{eq:kappaCond_app} is true,
then we can define a density $g_{\kappa_{\theta, g}}(y) = \left(\frac{f_\theta(y)}{f_{\theta_0}(y)}\right)^{\kappa_{\theta, g}} g(y)$,
and the test $\nu(f_\theta,f_{\theta_0})$ is equivalent to $\nu(g_{\kappa_{\theta, g}},g)$ with threshold $\kappa_{\theta, g} A$.
The estimate \eqref{eq:UB_SPRT_mis} is then just the classical estimate of the probability for a one-sided test to stop in finite time
under null hypothesis,
obtained by applying Theorem 1.1 in \cite{wood-nonlin-ren-th-book-1982} (also see Theorem 3.1 in \cite{wood-nonlin-ren-th-book-1982}).

We will now use \eqref{eq:UB_SPRT_mis} to obtain an upper bound on $\Prob_g(N < \infty)$.
Towards that end, we revisit Section 3 of \cite{lord-amstat-1971} and modify the proof there appropriately to suit out needs.
Note that
\[
\log \frac{f_\theta(y)}{f_{\theta_0}(y)} = (\theta-\theta_0) T(y) - b(\theta)+b(\theta_0).
\]
Thus, with
\[
S_m := \sum_{i=1}^m T(Y_i),
\]
we have
\begin{equation*}
\begin{split}
\sup_{\theta: |\theta-\theta_0| \geq \epsilon} \; &\sum_{i=1}^m \log \frac{f_\theta(Y_i)}{f_{\theta_0}(Y_i)}  \\
&= \sup_{\theta: |\theta-\theta_0| \geq \epsilon} \; (\theta-\theta_0) S_m - m(b(\theta)-b(\theta_0)).
\end{split}
\end{equation*}
Now,
\begin{equation}\label{eq:GLR_union}
\begin{split}
&\left\{\sup_{\theta: |\theta-\theta_0| \geq \epsilon} \;  (\theta-\theta_0) S_m - m(b(\theta)-b(\theta_0)) > A\right\} \\
&\quad=\left\{ S_m > \inf_{\theta: \theta \geq \theta_0 + \epsilon} \frac{A + m(b(\theta)-b(\theta_0))}{\theta-\theta_0} \right\}\\
&\quad \quad \cup \left\{ S_m < \sup_{\theta: \theta \leq \theta_0 - \epsilon} \frac{A + m(b(\theta)-b(\theta_0))}{\theta-\theta_0} \right\}.
\end{split}
\end{equation}
This is because if the left hand side is true, then there is $\theta_1$ such that $\{(\theta_1-\theta_0) S_m - m(b(\theta_1)-b(\theta_0)) > A\}$, and $\theta_1$ could be either greater or less than $\theta_0$, making $\theta_1-\theta_0$ positive or negative. Thus, left hand side is a subset of the right hand side. An identical argument given in reverse justifies that the right hand side is a subset of the left.

By Assumption~\ref{assum:kappa_gtheta}
there exists a positive constant $\kappa_{\theta, g}$ satisfying \eqref{eq:kappaCond} for every $\theta \in \Theta$
with $|\theta-\theta_0| \geq \epsilon$. Furthermore, there exists $\kappa_g$ such that
\begin{equation}\label{eq:kappastar_app}
0 < \kappa_g \leq \inf_{\theta \in \Theta: |\theta-\theta_0| \geq \epsilon} (\kappa_{\theta, g} ).
\end{equation}
With this assumption we have an upper bound on the estimate in \eqref{eq:UB_SPRT_mis}: $\forall \theta \in \Theta$, $|\theta-\theta_0| \geq \epsilon$,
\begin{equation}\label{eq:UB2_SPRT_mis}
\Prob_g(\nu(f_\theta,f_{\theta_0}) < \infty) \leq e^{-\kappa_{\theta, g} A} \leq e^{-\kappa_g A}.
\end{equation}

Now consider the infimum on the right hand side of \eqref{eq:GLR_union}. Let the infimum be approached along the sequence $\{\theta_\ell\}$.
Then,
\begin{equation}
\begin{split}
\Prob_g & \left\{ S_m > \inf_{\theta: \theta \geq \theta_0 + \epsilon} \frac{A + m(b(\theta)-b(\theta_0))}{\theta-\theta_0} \right\}\\
&=\lim_{\ell \to \infty} \Prob_g\left\{ S_m > \frac{A + m(b(\theta_\ell)-b(\theta_0))}{\theta_\ell-\theta_0} \right\}\\
&\leq\limsup_{\ell \to \infty} \Prob_g\left( \nu(f_{\theta_\ell}, f_{\theta_0}) \leq m\right)\\
&\leq\limsup_{\ell \to \infty} \Prob_g\left( \nu(f_{\theta_\ell}, f_{\theta_0}) < \infty\right)\\
&\leq e^{-\kappa_g A},
\end{split}
\end{equation}
where the last inequality follows from \eqref{eq:UB2_SPRT_mis}.
An almost identical argument yields the same bound on the probability of the other event involving a supremum in \eqref{eq:GLR_union}. Thus,
\begin{equation}
\begin{split}
\Prob_g & (N=m) \\
&\leq \Prob_g\left\{\sup_{\theta: |\theta-\theta_0| \geq \epsilon} \; (\theta-\theta_0) S_m - m(b(\theta)-b(\theta_0)) > A\right\} \\
&\leq 2 e^{-\kappa_g A}.
\end{split}
\end{equation}

By Assumption~\ref{assump:I_theta}, $I(\theta)$, the Kullback-Leibler divergence between $f_\theta$ and $f_{\theta_0}$,
increases with $|\theta-\theta_0|$.
Because of this assumption, if
\begin{equation}\label{eq:Abym_smaller}
\frac{A}{m} < \min\{I(\theta_0+\epsilon), I(\theta_0-\epsilon)\},
\end{equation}
then the infimum and supremum on the right hand side of \eqref{eq:GLR_union} are achieved at the boundaries $\theta_0+\epsilon$, and $\theta_0-\epsilon$, respectively.
To see this, we differentiate to show that the term inside the infimum is equal to
\begin{equation}
\begin{split}
\frac{d}{d\theta}&\frac{A + m(b(\theta)-b(\theta_0))}{\theta-\theta_0} \\
&\quad= \frac{m [(\theta-\theta_0) b'(\theta) - (b(\theta)-b(\theta_0))]-A}{(\theta-\theta_0)^2}\\
&\quad =\frac{m I(\theta)-A}{(\theta-\theta_0)^2}.
\end{split}
\end{equation}
Thus, setting the derivative to zero shows that the local interior minima $\theta^*$ must satisfy
\begin{equation}
I(\theta^*)=A/m.
\end{equation}
Since, $\Theta$ is assumed
to be an interval and the term inside the infimum is continuous, it must achieve its minimum on $[\theta_0+\epsilon, \theta_M]$,
where $\theta_M$ is the rightmost point of $\Theta$. The condition \eqref{eq:Abym_smaller} guarantees that the minimum cannot be achieved on $(\theta_0+\epsilon, \theta_M)$. Furthermore, it cannot be achieved at $\theta_M$ since otherwise
we would have the contradiction
\[
\frac{A}{m} < I(\theta_0 + \epsilon) \leq I(\theta_M) \leq \frac{A}{m},
\]
where the last inequality follows from the standard necessary condition on optimization over convex sets;
see Proposition 2.1.2 in \cite{bert_nlp_book_2004}. Almost identical arguments allows us to prove
that the supremum on the right hand side of \eqref{eq:GLR_union} is achieved at $\theta_0-\epsilon$ if
\eqref{eq:Abym_smaller} is true.

Define
$$
\mathcal{M}=\left\{m: m > \frac{A}{\min\{I(\theta_0+\epsilon), I(\theta_0-\epsilon)\}}\right\}.
$$
We have the estimate
\begin{equation}
\begin{split}
&\Prob_g \left(\frac{A}{\min\{I(\theta_0+\epsilon), I(\theta_0-\epsilon)\}} < N < \infty\right) \\
&=\Prob_g \left(\cup_{m \in \mathcal{M}} \{N=m\}\right) \\
\leq &\Prob_g\left\{\bigcup_{m \in \mathcal{M}} \left\{\sup_{\theta: |\theta-\theta_0| \geq \epsilon} \hspace{-0.3cm}(\theta-\theta_0) S_m - m(b(\theta)-b(\theta_0)) > A\right\}\right\} \\
\leq &\Prob_g\left\{\bigcup_{m \in \mathcal{M}} \left\{S_m > \inf_{\theta: \theta \geq \theta_0 + \epsilon} \frac{A + m(b(\theta)-b(\theta_0))}{\theta-\theta_0} \right\}\right\} \\
&+\Prob_g\left\{\bigcup_{m \in \mathcal{M}} \left\{S_m < \sup_{\theta: \theta \leq \theta_0 - \epsilon} \frac{A + m(b(\theta)-b(\theta_0))}{\theta-\theta_0} \right\}\right\} \\
= &\Prob_g\left\{\bigcup_{m \in \mathcal{M}} \left\{S_m > \frac{A + m(b(\theta_0+\epsilon)-b(\theta_0))}{\theta_0+\epsilon-\theta_0} \right\}\right\} \\
&+\Prob_g\left\{\bigcup_{m \in \mathcal{M}} \left\{S_m < \frac{A + m(b(\theta_0-\epsilon)-b(\theta_0))}{\theta_0-\epsilon-\theta_0} \right\}\right\} \\
\leq &\Prob_g\left\{\bigcup_{m \in \mathcal{M}} \left\{\nu(f_{\theta_0+\epsilon}, f_{\theta_0}) \leq m \right\}\right\} \\
&+\Prob_g\left\{\bigcup_{m \in \mathcal{M}} \left\{\nu(f_{\theta_0-\epsilon}, f_{\theta_0})\leq m \right\}\right\} \\
 \leq & \Prob_g\left(\nu(f_{\theta_0+\epsilon}, f_{\theta_0}) < \infty\right) +
\Prob_g\left(\nu(f_{\theta_0-\epsilon}, f_{\theta_0}) < \infty\right)\\
\leq & 2 e^{-\kappa_g A}.
\end{split}
\end{equation}

Thus, similar to the estimate in \cite{lord-amstat-1971}, we have the estimate
\begin{equation}
\begin{split}
\Prob_g & (N<\infty) \\
&=\sum_{m=1}^{\lfloor\frac{A}{\min\{I(\theta_0+\epsilon), I(\theta_0-\epsilon)\} } \rfloor} \Prob_g (N=m) \\
&\quad +\Prob_g \left(\frac{A}{\min\{I(\theta_0+\epsilon), I(\theta_0-\epsilon)\}} < N < \infty\right)\\
&\leq 2 e^{-\kappa_g A} \left(\frac{A}{\min\{I(\theta_0+\epsilon), I(\theta_0-\epsilon)\}} + 1\right).
\end{split}
\end{equation}
From \eqref{eq:onesided_toCuSum}
\begin{equation}
\begin{split}
\mathsf{E}_{g}[\taug]\geq  \frac{1}{2 e^{-\kappa_g A} \left(\frac{A}{\min\{I(\theta_0+\epsilon), I(\theta_0-\epsilon)\}} + 1\right)}.
\end{split}
\end{equation}
This proves the first part of the theorem. The second part is now obvious.
\end{IEEEproof}

\begin{IEEEproof}[Proof of Theorem~\ref{thm:Misspec_Delay}]
Let $\theta_g$ be as in Assumption~\ref{assump:delay}. Then note that
\begin{equation}
\begin{split}
\taug &= \inf\left\{m\geq 1: \max_{1\leq k \leq m} \; \sup_{\theta: |\theta-\theta_0| \geq \epsilon} \; \sum_{i=k}^m \log \frac{f_\theta(Y_i)}{f_{\theta_0}(Y_i)} > A\right\}\\
&\leq \inf\left\{m\geq 1: \sum_{i=1}^m \log \frac{f_{\theta_g}(Y_i)}{f_{\theta_0}(Y_i)} > A\right\}.
\end{split}
\end{equation}
Assumption~\ref{assump:delay} implies that the drift of the random walk with increments
$\log \frac{f_{\theta_g}(Y_i)}{f_{\theta_0}(Y_i)}$ is positive when samples are drawn from $g$. The theorem now follows from Proposition 8.21 in \cite{sieg-seq-anal-book-1985}: as $A \to \infty$
\begin{equation}
\begin{split}
\Expect_g[\taug] \leq \frac{A}{\Expect_g\left( \log \frac{f_{\theta_g}(Y_1)}{f_{\theta_0}(Y_1)}\right)}(1+o(1)).
\end{split}
\end{equation}
\end{IEEEproof}

\begin{IEEEproof}[Proof of Lemma~\ref{lem:Gaussian}]
For the  case of a Gaussian distribution \eqref{eq:GaussDenAssumptions} with mean parameter $\theta$ the KL divergence KL$(f_\theta\|f_{\theta_0})$ and  the integral \eqref{eq:kappaCond} have closed form analytical expressions. The KL divergence can be shown to be quadratic function of $\theta$ taking its minimum at $\theta_0$, so Assumption 3 is satisfied. 
\[
\int \left(\frac{f_\theta(y)}{f_{\theta_0}(y)}\right)^{\kappa_{\theta, g}} f_{\tilde{\theta}_0}(y) \; dy= 1
\]
can be explicitly solved for $\kappa_{\theta, g}$ giving two solutions: $\kappa_{\theta, g}=0$ and
that given by \eqref{eq:gauss_kap_g_th}.
The latter is positive only if $|\tilde{\theta}_0-\theta_0|< \epsilon/2$. This proves the first part of the theorem.

The second part is true because $\kappa_{\theta, g}$ is monotonic in $\theta$, and its value is smallest when
$\theta$ is either equal to $\theta_0+\epsilon$ or  $\theta_0-\epsilon$. This value has the explicit expression given by
the right most expression in \eqref{eq:gauss_kap_g}.

The third part of the theorem is true because the expression for $\kappa_g$ given in \eqref{eq:gauss_kap_g} is monotonic in
$|\tilde{\theta}_0-\theta_0|$.
\end{IEEEproof}

\balance


\bibliographystyle{ieeetr}



\bibliography{QCD_verSubmitted}

\begin{thebibliography}{10}

\bibitem{hero-bala-submitted-2015}
A.~Hero and B.~Rajaratnam, ``Foundational principles for large scale inference:
  Illustrations through correlation mining,'' {\em IEEE Proceedings}, 2015.
\newblock In press. \url{http://arxiv.org/abs/1505.02475}.

\bibitem{veer-bane-elsevierbook-2013}
V.~V. Veeravalli and T.~Banerjee, {\em Quickest Change Detection}.
\newblock Academic Press Library in Signal Processing: Volume 3 -- Array and
  Statistical Signal Processing, 2014.
\newblock \url{http://arxiv.org/abs/1210.5552}.

\bibitem{poor-hadj-qcd-book-2009}
H.~V. Poor and O.~Hadjiliadis, {\em Quickest detection}.
\newblock Cambridge University Press, 2009.

\bibitem{tart-niki-bass-2014}
A.~G. Tartakovsky, I.~V. Nikiforov, and M.~Basseville, {\em Sequential
  Analysis: {Hypothesis} Testing and Change-Point Detection}.
\newblock Statistics, CRC Press, 2014.

\bibitem{wald_book_1947}
A.~Wald, {\em Sequential analysis}.
\newblock Dover Publication, 2013.

\bibitem{wald-wolf-amstat-1948}
A.~Wald and J.~Wolfowitz, ``Optimum character of the sequential probability
  ratio test,'' {\em Ann. Math. Statist.}, vol.~19, no.~3, pp.~pp. 326--339,
  1948.

\bibitem{girs_rubin_1952}
M.~A. Girshick and H.~Rubin, ``A {B}ayes approach to a quality control model,''
  {\em Ann. Math. Statist.}, pp.~114--125, 1952.

\bibitem{shir-siamtpa-1963}
A.~N. Shiryaev, ``On optimum methods in quickest detection problems,'' {\em
  Theory of Prob and App.}, vol.~8, pp.~22--46, 1963.

\bibitem{tart-veer-siamtpa-2005}
A.~G. Tartakovsky and V.~V. Veeravalli, ``General asymptotic {Bayesian} theory
  of quickest change detection,'' {\em SIAM Theory of Prob. and App.}, vol.~49,
  pp.~458--497, Sept. 2005.

\bibitem{page-biometrica-1954}
E.~S. Page, ``Continuous inspection schemes,'' {\em Biometrika}, vol.~41,
  pp.~100--115, June 1954.

\bibitem{lord-amstat-1971}
G.~Lorden, ``Procedures for reacting to a change in distribution,'' {\em Ann.
  Math. Statist.}, vol.~42, pp.~1897--1908, Dec. 1971.

\bibitem{mous-astat-1986}
G.~V. Moustakides, ``Optimal stopping times for detecting changes in
  distributions,'' {\em Ann. Statist.}, vol.~14, pp.~1379--1387, Dec. 1986.

\bibitem{ritov-astat-1990}
Y.~Ritov, ``Decision theoretic optimality of the {CUSUM} procedure,'' {\em Ann.
  Statist.}, vol.~18, pp.~1464--1469, Nov. 1990.

\bibitem{lai-ieeetit-1998}
T.~L. Lai, ``Information bounds and quick detection of parameter changes in
  stochastic systems,'' {\em IEEE Trans. Inf. Theory}, vol.~44, pp.~2917
  --2929, Nov. 1998.

\bibitem{lord-poll-nonantiest-2005}
G.~Lorden and M.~Pollak, ``Nonanticipating estimation applied to sequential
  analysis and changepoint detection,'' {\em Ann. Statist.}, pp.~1422--1454,
  2005.

\bibitem{poll-astat-1985}
M.~Pollak, ``Optimal detection of a change in distribution,'' {\em Ann.
  Statist.}, vol.~13, pp.~206--227, Mar. 1985.

\bibitem{hero-bala-IT-2012}
A.~Hero and B.~Rajaratnam, ``Hub discovery in partial correlation graphs,''
  {\em IEEE Trans. Inf. Theory}, vol.~58, no.~9, pp.~6064--6078, 2012.

\bibitem{lai-jrss-1995}
T.~L. Lai, ``Sequential changepoint detection in quality control and dynamical
  systems,'' {\em J. Roy. Statist. Soc. Suppl.}, vol.~57, no.~4, pp.~pp.
  613--658, 1995.

\bibitem{ander-mutlstat-1996}
T.~W. Anderson, {\em An Introduction to Multivariate Statistical Analysis}.
\newblock New York, NY: Wiley, 2003.

\bibitem{AndersonTR}
T.~W. Anderson, ``Non-normal multivariate distributions: inference based on
  elliptically contoured distributions,'' tech. rep., July 1992.

\bibitem{hero-bala-jasta-2011}
A.~Hero and B.~Rajaratnam, ``Large-scale correlation screening,'' {\em J. Amer.
  Statist. Assoc.}, vol.~106, no.~496, pp.~1540--1552, 2011.

\bibitem{moran1980testing}
P.~Moran, ``Testing the largest of a set of correlation coefficients,'' {\em
  Australian \& New Zealand Journal of Statistics}, vol.~22, no.~3,
  pp.~289--297, 1980.

\bibitem{cameron1985new}
M.~Cameron and G.~Eagleson, ``A new procedure for assessing large sets of
  correlations,'' {\em Australian \& New Zealand Journal of Statistics},
  vol.~27, no.~1, pp.~84--95, 1985.

\bibitem{anandkumar2009detection}
A.~Anandkumar, L.~Tong, and A.~Swami, ``Detection of gauss--markov random
  fields with nearest-neighbor dependency,'' {\em IEEE Transactions on
  Information Theory}, vol.~55, no.~2, pp.~816--827, 2009.

\bibitem{cox2014multivariate}
D.~R. Cox and N.~Wermuth, {\em Multivariate dependencies: Models, analysis and
  interpretation}.
\newblock Chapman and Hall/CRC, 2014.

\bibitem{robe-technometrics-1966}
S.~W. Roberts, ``A comparison of some control chart procedures,'' {\em
  Technometrics}, vol.~8, pp.~411--430, Aug. 1966.

\bibitem{tart-etal-thirdord-2012}
A.~G. Tartakovsky, M.~Pollak, and A.~S. Polunchenko, ``Third-order asymptotic
  optimality of the generalized {Shiryaev-Roberts} changepoint detection
  procedures,'' {\em Theory of Prob and App.}, vol.~56, no.~3, pp.~457--484,
  2012.

\bibitem{wood-nonlin-ren-th-book-1982}
M.~Woodroofe, {\em Nonlinear Renewal Theory in Sequential Analysis}.
\newblock CBMS-NSF regional conference series in applied mathematics, SIAM,
  1982.

\bibitem{sricharan2011local}
K.~Sricharan, A.~O. Hero, and B.~Rajaratnam, ``A local dependence measure and
  its application to screening for high correlations in large data sets,'' in
  {\em Information Fusion (FUSION), 2011 Proceedings of the 14th International
  Conference on}, pp.~1--8, IEEE, 2011.

\bibitem{bert_nlp_book_2004}
D.~P. Bertsekas, {\em Nonlinear programming}.
\newblock Belmont, MA: Athena Scientific, 2004.

\bibitem{sieg-seq-anal-book-1985}
D.~Siegmund, {\em Sequential Analysis: {Tests} and Confidence Intervals}.
\newblock Springer series in statistics, Springer-Verlag, 1985.

\end{thebibliography}

\end{document}